\documentclass{siamltex}
\usepackage{a4wide,amssymb,graphics,graphicx,textcomp}
\usepackage{enumerate,textcomp,multirow,amsmath}
\newtheorem{remark}[theorem]{Remark}
\newtheorem{example}[theorem]{Example}

\numberwithin{equation}{section}
\numberwithin{table}{section}
\numberwithin{figure}{section}
\newcommand{\eproof}{\space
    {\ \vbox{\hrule\hbox{\vrule height1.3ex\hskip0.8ex\vrule}\hrule}}\par}

\newcommand {\mat}  [1] {\left[\begin{array}{#1}}
\newcommand {\rix}      {\end{array}\right]}

\catcode`@=11     
\font\tenex=cmex10 
\newdimen\p@renwd
\setbox0=\hbox{\tenex B} \p@renwd=\wd0 
\def\bmat#1{\begingroup \m@th
  \setbox\z@\vbox{\def\cr{\crcr\noalign{\kern2\p@\global\let\cr\endline}}%
    \ialign{$##$\hfil\kern2\p@\kern\p@renwd&\thinspace\hfil$##$\hfil
      &&\quad\hfil$##$\hfil\crcr
      \omit\strut\hfil\crcr\noalign{\kern-\baselineskip}%
      #1\crcr\omit\strut\cr}}%
  \setbox\tw@\vbox{\unvcopy\z@\global\setbox\@ne\lastbox}%
  \setbox\tw@\hbox{\unhbox\@ne\unskip\global\setbox\@ne\lastbox}%
  \setbox\tw@\hbox{$\kern\wd\@ne\kern-\p@renwd\left[\kern-\wd\@ne
    \global\setbox\@ne\vbox{\box\@ne\kern2\p@}%
    \vcenter{\kern-\ht\@ne\unvbox\z@\kern-\baselineskip}\,\right]$}%
  \null\;\vbox{\kern\ht\@ne\box\tw@}\endgroup}
\catcode`@=12    

\newcommand{\wt}{\widetilde}

\newcommand{\backmatter}{\appendix
\def\chaptermark##1{\markboth{%
\ifnum  \c@secnumdepth > \m@ne  \@chapapp\ \thechapter:  \fi  ##1}{%
\ifnum  \c@secnumdepth > \m@ne  \@chapapp\ \thechapter:  \fi  ##1}}%
\def\sectionmark##1{\relax}}

\makeatletter
\newcommand*{\rom}[1]{\expandafter\@slowromancap\romannumeral #1@}
\makeatother

\newif\ifMDlatex
\catcode`@=11     
\def\MD@us#1{\csname#1style\endcsname}
\def\MD@uf#1{\csname#1font\endcsname}
\def\MD@t{text}
\def\MD@s{script}
\def\MD@ss{scriptscript}
\newdimen\MD@unit
\MD@unit\p@
\def\MD@changestyle#1{
  \relax\MD@unit0.1\fontdimen6\MD@uf{#1}0
  \everymath\expandafter{\the\everymath\MD@us{#1}}
}
\def\MD@dot{$\m@th\ldotp$}
\def\MD@palette#1{\mathchoice{#1\MD@t}{#1\MD@t}{#1\MD@s}{#1\MD@ss}}
\def\MD@ddots#1{{\MD@changestyle{#1}%
  \mkern1mu\raise7\MD@unit\vbox{\kern7\MD@unit\hbox{\MD@dot}}%
  \mkern2mu\raise4\MD@unit\hbox{\MD@dot}%
  \mkern2mu\raise \MD@unit\hbox{\MD@dot}\mkern1mu}}%
\def\MD@iddots#1{{\MD@changestyle{#1}%
  \mkern1mu\raise \MD@unit\hbox{\MD@dot}%
  \mkern2mu\raise4\MD@unit\hbox{\MD@dot}%
  \mkern2mu\raise7\MD@unit\vbox{\kern7\MD@unit\hbox{\MD@dot}}}}%
\def\MD@vdots#1{\vbox{\MD@changestyle{#1}%
    \baselineskip4\MD@unit\lineskiplimit\z@
    \kern6\MD@unit\hbox{\MD@dot}\hbox{\MD@dot}\hbox{\MD@dot}}}%
\ifMDlatex
  \DeclareRobustCommand\ddots{\mathinner{\MD@palette\MD@ddots}}%
  \DeclareRobustCommand\iddots{\mathinner{\MD@palette\MD@iddots}}%
  \DeclareRobustCommand\vdots{\mathinner{\MD@palette\MD@vdots}}%
\else
  \def\ddots{\mathinner{\MD@palette\MD@ddots}}%
  \def\iddots{\mathinner{\MD@palette\MD@iddots}}%
  \def\vdots{\mathinner{\MD@palette\MD@vdots}}%
\fi
\catcode`@=12    

\newcommand {\comment}[1]{} 
\newcommand{\C}{{\mathbb C}}
\newcommand{\R}{{\mathbb R}}

\begin{document}
\title{Structured eigenvalue/eigenvector backward errors of matrix pencils arising in optimal control}

\author{Christian Mehl\thanks{Institut f${\rm \ddot{u}}$r Mathematik, MA 4-5 TU Berlin, Str.\@ d.\@ 17.\@ Juni 136, D-10623 Berlin, Germany. Email: \{mehrmann, mehl\}@math.tu-berlin.de.  C. Mehl and V. Mehrmann gratefully acknowledge support from Einstein Center ECMath via project SE3: Stability analysis of power networks and power network models. V. Mehrmann also acknowledges support Deutsche Forschungsgemeinschaft through CRC 910 \emph{Control of Self-Organizing Nonlinear Systems} via project A02: Analysis and computation of stability exponents for delay differential-algebraic equations.}
\and Volker Mehrmann$^{*}$
\and Punit Sharma\thanks{Department of Mathematics and Operational Research, University of Mons, Rue de Houdain 9, 7000 Mons, Belgium. Email: punit.sharma@umons.ac.be. Punit Sharma acknowledges the support of the ERC starting grant n$^\text{o}$ 679515.}
}

\maketitle

\begin{abstract}
Eigenvalue and eigenpair backward errors are computed for matrix pencils arising in optimal control.
In particular, formulas for backward errors are developed that are obtained under block-structure-preserving
and symmetry-structure-preserving perturbations. It is shown that these eigenvalue and eigenpair
backward errors are sometimes significantly larger than the corresponding backward errors that are
obtained under perturbations that ignore the special structure of the pencil.
\end{abstract}


\section{Introduction}
In this paper we consider the perturbation theory, in particular the calculation of structured backward errors,  for eigenvalues and eigenvectors of structured matrix pencils $L(z)$ of the form
\begin{equation}\label{eq:defL_cz2}
L(z)=M+zN:=\mat{ccc} 0 & J-R & B\\ (J-R)^H & Q& 0 \\B^H& 0 &S\rix +z \mat{ccc}0&E&0\\-E^H &0&0\\0&0&0\rix
\end{equation}
where $J,R,E,Q \in \C^{n, n}$, $B \in \C^{n,m} $ and $S\in \C^{m, m}$ satisfy
$J^H=-J$, $R^H=R$
$E^H=E$, $Q^H=Q$, and $S^H=S > 0$, i.e., $S$ is positive definite.
These pencils are special cases of so-called \emph{even pencils}, i.e.,
matrix pencils $P(\lambda)$ satisfying $P(z=P(-z)^H$, see, e.g., \cite{MacMMM06b}.
Even pencils with an additional block-structure as in~\eqref{eq:defL_cz2}
arise in optimal control and $H_\infty$ control problems as well as in the passivity analysis of dynamical systems. For instance, if one considers the optimal control problem of minimizing the cost functional
\[
\int_{t_0}^{\infty} x^H Q  x + u^H S u \ dt
\]
subject to the constraint
\begin{equation}\label{eq:system}
E\dot x =A x+Bu,\quad x(t_0)=x^0,
\end{equation}
then it is well known, see \cite{KunM08,Meh91}, that the optimal solution is associated with the deflating subspace of a pencil of the form~\eqref{eq:defL_cz2} associated with the finite eigenvalues in the open left half plane. If there exist exactly $n$ eigenvalues in the open left half plane then
this deflating subspace is an \emph{extended Lagrangian subspace}. (For other applications in passivity analysis and robust control,  see \cite{BenLMV15}.)
Note that for general descriptor systems we need not have that $E=E^H$. However, if this is not the case then we can just carry out the polar decomposition, see \cite{GolV96}, to obtain $E= U \tilde E$ with $U$ unitary and $\tilde E =\tilde E^H$. Multiplying equation~\eqref{eq:system} from the left with $U^H$ we obtain a new system that has the desired property $E=E^H$, so w.l.o.g. we assume that $E=E^H$ and then partition in its skew-symmetric and symmetric part $A=J-R$.
Note that this condition automatically holds if
\eqref{eq:system} is  a \emph{port-Hamiltonian descriptor system}, see \cite{BeaMXZ17_ppt,Sch13}. In this case we furthermore have that $R\geq 0$, i.e. is positive semidefinite.

The solution of the optimal control problem becomes highly ill-conditioned when eigenvalues are close the imaginary axis and the solution usually seizes to exist when the eigenvalues are on the imaginary axis \cite{BenBMX07,FreMX02}. When eigenvalues on the imaginary axis exist then it is an important question to find small perturbations to the system~(\ref{eq:system}) or the pencil~(\ref{eq:defL_cz2}) that remove the eigenvalues from the imaginary axis \cite{AlaBKMM11,GilMS17}. These questions motivate  the principle aims of this paper to determine  backward errors
associated with eigenvalues on the imaginary axis of pencils of the form~\eqref{eq:defL_cz2}.
We will consider in this paper the special case of  pencils with $Q=0$, which arises in optimal control without state weighting,  and in the context of passivity analysis~\cite{GenHNSVX02,GilMS17}. Thus, we will
consider a pencil of the form
\begin{equation}\label{eq:defL_cz}
L(z)=M+zN:=\mat{ccc} 0 & J-R & B\\ (J-R)^H & 0& 0 \\B^H& 0 &S\rix +z \mat{ccc}0&E&0\\-E^H &0&0\\0&0&0\rix
\end{equation}

In the following, $\|\cdot\|$ denotes the spectral norm of a vector or a matrix and ${\|A\|}_{F}$ denotes the Frobenius norm of a matrix $A$.  ${\rm Herm}(n)$ and ${\rm SHerm}(n)$ respectively denote the set of Hermitian and skew-Hermitian matrices of size $n$. By $i\R$ we denote the set of nonzero purely imaginary numbers, i.e.,   $i\R=\left\{i\alpha~|~\alpha \in \R\setminus\{0\} \right\}$, and by $I_n$ the identity matrix of size $n$. For a matrix $A$ we write $A=0$ if each entry of $A$ is equal to zero.

The sensitivity analysis of eigenvalues and eigenvalue/eigenvector pairs (in the following called eigenpairs) of matrix pencils and matrix polynomials
with various structures has recently received  a lot of attention, see, e.g., \cite{Adh08,AhmA09,AdhA09,HigH98,KarKT06,Tis03}. In particular, backward error formulas for structured matrix pencils and  polynomials with respect to structure preserving perturbations have been obtained in~\cite{Adh08,AdhA09} and in~\cite{BorKMS14,BorKMS15}, respectively.

For pencils of the form~\eqref{eq:defL_cz}, if the structure of the pencil is ignored, then for a given pair
$(\lambda,x)\in \mathbb C\times (\C^{2n+m}\setminus\{0\})$ the \emph{eigenpair backward error} is defined as
\begin{eqnarray*}
\eta(L,\lambda,x)=\inf\Big\{\|[\Delta_M~\,\Delta_N]\|_{F}\, \Big |\,
\Delta_M,\, \Delta_N \in \C^{2n+m , 2n+m},\big((M-\Delta_M)+\lambda(N-\Delta_N)\big)x=0\Big\}.
\end{eqnarray*}
It can be interpreted as the Frobenius norm of the smallest perturbation that makes
$(\lambda,x)$ being an eigenpair of the perturbed pencil. Minimizing this expression over all
$(\lambda,x)\in (i\mathbb R)\times(\mathbb C^{2n+m})$ we obtain the distance of $L(z)$
to the next pencil having eigenvalues on the imaginary axis and thus, the passivity radius of $L(z)$,
see \cite{GilS17b,OveV05}. If the even structure of the pencil
is taken into account, then a structured eigenpair backward error with respect to
structure-preserving perturbations can be defined as
\begin{eqnarray*}
\eta^{even}(L,\lambda,x)=\inf\Big\{\|[\Delta_M~\,\Delta_N]\|_{F} \,\Big |~
\Delta_M\in {\rm Herm}(2n+m),\,\Delta_N \in {\rm SHerm}(2n+m)\quad\\
\big((M-\Delta_M)+\lambda(N-\Delta_N)\big)x=0\Big\}.
\end{eqnarray*}
Clearly, we have $\eta(L,\lambda,x) \leq \eta^{even}(L,\lambda,x)$. In fact, for a given
$(\lambda,x)\in \C \times  (\C^{2n+m}\setminus \{0\})$, it is well known by \cite[Theorem 4.6]{AhmA09} that
\begin{equation}\label{17.9.17one}
\eta(L,\lambda,x)=\frac{\|L(\lambda)x\|}{\|x\|\sqrt{1+|\lambda|^2}},
\end{equation}
and by \cite[Theorem 3.3.7]{Adh08} that
\begin{equation}\label{17.9.17two}
\eta^{even}(L,\lambda,x)=\sqrt{\frac{2{\|x\|}^2{\|L(\lambda)x\|}^2 - {|x^HL(\lambda)x|}^2}
{{\|x\|}^4 (1+|\lambda|^2)}}.
\end{equation}
However, both formulas ignore the special block-structure of the pencil $L(z)$, in particular the zero structure and the definiteness of the matrix $S$, and as we will show
in this paper, eigenpair backward errors with respect to perturbations that preserve the
block-structure and possibly also the symmetry-structure may be significantly larger
than the more generally obtained backward errors $\eta(L,\lambda,x)$ and $\eta^{even}(L,\lambda,x)$.

The remainder of this paper is organized as follows. In section 2, we review some minimal norm mapping problems.
In section 3, we introduce a terminology and define block- and symmetry-structure-preserving eigenpair or
eigenvalue backward errors for pencils $L(z)$ of the form~\eqref{eq:defL_cz}. These backward errors are
computed while perturbing any two, three or all of the blocks $J,R,E$ or $B$ in sections 4, 5 and 6, respectively.
The significance of these block- and symmetry-structure-preserving backward errors over $\eta(L,\lambda,x)$ and $\eta^{even}(L,\lambda,x)$
is shown via some numerical examples in section 7.

\section{Preliminaries}

An important tool for the computation of backward errors are minimal norm solutions
to mapping problems. In this section we will review some of these results  and restate them in a form that we need in the
following sections.

The solution to the \emph{skew-Hermitian mapping problem} to find
$\Delta \in {\rm SHerm}(n)$ that maps a matrix $X\in \C^{n , k}$ to $Y\in \C^{n , k} $ is well known, see, e.g., \cite{Adh08}, where also
solutions that are minimal with respect to the spectral and the Frobenius norm are characterized.
The following theorem is a particular case of \cite[Theorem 2.2.3]{Adh08}.
\begin{theorem}\label{eq:SHerm_map}
Let $X,\,Y \in \C^{n , k}$. Then there exist
$\Delta \in {\rm SHerm}(n)$ satisfying $\Delta X=Y$ if and only if $YX^{\dag}X=Y$
and $Y^HX=-X^HY$. If the latter conditions are satisfied then
\begin{equation*}
\min\big\{\|\Delta\|_F\,\big|\,\Delta \in {\rm SHerm}(n),\,\Delta X =Y\big\}=\sqrt{2 {\|YX^\dag\|}_F^2-
\operatorname{trace}\big(YX^\dag(YX^\dag)^H(XX^\dag)\big)}
\end{equation*}
 and the unique minimum is attained for
\begin{equation*}
\widehat \Delta= YX^\dag-(YX^\dag)^H-(X^\dag)^HX^HYX^\dag.
\end{equation*}
\end{theorem}
The second mapping problem that we will need is the following, see \cite{KahPJ82}.
%
\begin{theorem} \label{thm:spcl_matrix_map}
Let $u\in \C^m$, $r\in \C^{n}\setminus \{0\}$, $w\in \C^n$ and $s\in \C^m\setminus \{0\}$. Define
\[
\mathcal S=\{\Delta \in \C^{n,m}~|~\Delta u=r,~\Delta^Hw=s\}.
\]
Then $\mathcal S \neq \emptyset$ if and only if $u^Hs=r^Hw$. If the latter condition is satisfied then
 \[
 \widehat \Delta=\frac{ru^H}{{\|u\|}^2}+\frac{ws^H}{{\|w\|}^2}-\frac{(s^Hu)wu^H}{{\|w\|}^2{\|u\|}^2}
\]
is the unique matrix such that $\widehat \Delta u=r$ and $\widehat \Delta^H w=s$, and
\[
\inf_{\Delta \in \mathcal S} {\|\Delta\|}_F={\|\widehat \Delta\|}_F=
\sqrt{\frac{{\|r\|}^2}{{\|u\|}^2}+\frac{{\|s\|}^2}{{\|w\|}^2}-\frac{{|s^Hu|}^2}{\|w\|\|u\|}}.
\]
Moreover,
\[
\inf_{\Delta \in \mathcal S}\|\Delta\|=\max\left\{\frac{\|r\|}{\|u\|},\frac{\|s\|}{\|w\|}\right\}.
\]
\end{theorem}
The following result, see \cite[Remark 2.4]{MehMS17} gives a real minimal Frobenius norm solution of
the mapping problem considered in Theorem~\ref{thm:spcl_matrix_map}.
\begin{theorem}\label{cor:real_spcl_cor}
Let $u\in \C^m$, $r\in \C^{n}$, $w\in \C^n$ and $s\in \C^m$ such that
$\operatorname{rank}([u~\bar u])=2$ and $\operatorname{rank}([w~\bar w])=2$ and define
\[
\mathcal S_\R=\{\Delta \in \R^{n,m}~|~\Delta u=r,~\Delta^Hw=s\}.
\]
Then $\mathcal S_\R \neq \emptyset$ if and only if $u^Hs=r^Hw$ and $u^Ts=r^Tw$. If the latter conditions
are satisfied, then
\[
\inf_{\Delta \in \mathcal S_\R}\|\Delta\|_F=\|\widetilde\Delta\|,
\]
where
 \[
 \widetilde \Delta=[r~\bar r][u~\bar u]^\dagger+([s~\bar s][w~\bar w]^\dagger)^H-
 ([s~\bar s][w~\bar w]^\dagger)^H[u~\bar u][u~\bar u]^\dagger.
 \]
\end{theorem}

We mention that the form of the minimal norm perturbation given in \cite[Remark 2.4]{MehMS17}
slightly differs from the one given here, because in \cite{MehMS17} it was presented using
real and imaginary parts rather than complex vectors and their complex conjugates.


\section{Structured eigenpair backward errors}\label{sec:section_def_epbr}
In this section we consider structured matrix pencils $L(z)$ of the form~\eqref{eq:defL_cz}.
We use the results on the mapping problems from the previous section to estimate structure-preserving
backward errors for eigenvalues $\lambda$ or eigenpairs $(\lambda,x)$ of $L(z)$, while perturbing only certain block entries
of $L(z)$ for the case when $\lambda$ is purely imaginary and $S$
is definite.
To distinguish between different cases, we introduce a terminology for perturbations $\Delta_M+z\Delta_N$ of the pencil $L(z)=M+zN$ that affect only some of the blocks $J,R,E,B$ of $L(z)$.
For example, suppose that only the blocks $J$ and $E$ in $L(z)$ are subject to perturbations. Then the corresponding
perturbations to $M$ and $N$ are given by
\begin{equation}\label{eq:M_N_def}
\Delta_M=\mat{ccc}0 & \Delta_J& 0\\ \Delta_J^H &0&0 \\ 0&0& 0\rix
\quad \text{and} \quad
\Delta_N=\mat{ccc}0 & \Delta_E& 0\\ -\Delta_E^H &0&0 \\ 0&0& 0\rix,
\end{equation}
where $\Delta_J,\, \Delta_E \in \C^{n,n}$. For $\lambda \in \C$ and $x \in \C^{2n+m}\setminus \{0\}$ we then define
\begin{enumerate}
\item the \emph{block-structure-preserving eigenpair backward error} $\eta^{\mathcal B}(J,E,\lambda,x)$
with respect to perturbations only to $J$ and $E$ by
\begin{equation}\label{eq:def_JE_eperr}
\eta^{\mathcal B}(J,E,\lambda,x)=\inf\Big\{\|[\Delta_J~\Delta_E]\|_{F} \,\Big |~
\big((M-\Delta_M)+\lambda(N-\Delta_N)\big)x=0,\,\Delta_M + z \Delta_N \in \mathcal B\Big\},
\end{equation}
where $\mathcal B$ denotes the set of all pencils $\Delta_M+z\Delta_N$ as in~\eqref{eq:M_N_def}
with $\Delta_J,\, \Delta_E \in \C^{n,n}$;
\item the \emph{symmetry-structure-preserving eigenpair backward error} $\eta^{\mathcal S}(J,E,\lambda,x)$
with respect to structure-preserving perturbations only to $J$ and $E$ by
\begin{equation}\label{eq:def_strJE_eperr}
\eta^{\mathcal S}(J,E,\lambda,x)=\inf\Big\{\|[\Delta_J~\Delta_E]\|_{F} \,\big |~
\big((M-\Delta_M)+\lambda(N-\Delta_N)\big)x=0,\,\Delta_M + z \Delta_N \in \mathcal S\Big\},
\end{equation}
where $\mathcal S$ denotes the set of all pencils $\Delta_M+z\Delta_N$ as in~\eqref{eq:M_N_def}
with $\Delta_J \in {\rm SHerm}(n)$ and $\Delta_E \in {\rm Herm}(n)$.
\end{enumerate}
For a given $\lambda \in \mathbb C$, we also define the \emph{block-structure-preserving} and \emph{symmetry-structure-preserving eigenvalue backward errors} $\eta^{\mathcal B}(J,E,\lambda)$ and
$\eta^{\mathcal S}(J,E,\lambda)$, respectively, by
\begin{equation*}
\eta^{\mathcal B}(J,E,\lambda):=\inf_{x \in \C^{2n+m}\setminus \{0\}}\eta^{\mathcal B}(J,E,\lambda,x)
\quad\text{and}\quad \eta^{\mathcal S}(J,E,\lambda):=\inf_{x \in \C^{2n+m}\setminus \{0\}}\eta^{\mathcal S}(J,E,\lambda,x).
\end{equation*}
For other combinations of perturbations to the blocks $J,R,E,B$ in $L(z)$,
the corresponding sets $\mathcal B$ and $\mathcal S$ as well as the block- and symmetry-structure-preserving
eigenpair or eigenvalue backward errors are defined analogously.

\section{Perturbation in any two of the blocks $J,R,E$, or $B$}

In this section, we compute block- and symmetry-structure-preserving  backward errors of
 $\lambda \in i\R$ and $ x \in \C^{2n+m}\setminus \{0\}$ as approximate eigenpair, resp. eigenvalue of the pencil $L(z)$ defined in~\eqref{eq:defL_cz} while perturbing any two of the blocks $J,R,E$, or $B$ at a time.

\subsection{Perturbation only in $J$ and $E$}
Let $L(z)$ be a pencil as in~\eqref{eq:defL_cz} and furthermore let $(\lambda,x) \in \C \times (\C^{2n+m}\setminus \{0\})$.
Suppose that only the blocks $J$ and $E$ of $L(z)$ are subject to perturbations.
Then by Section~\ref{sec:section_def_epbr}, $\mathcal B$ is the set of all pencils
$\Delta L(z)=\Delta_M + z \Delta_N$, where $\Delta_M$ and $\Delta_N$ have the block structure as in~\eqref{eq:M_N_def},
and $\mathcal S$ is the set of all pencils from $\mathcal B$ where in addition we have
$\Delta J^H=-\Delta J$ and $\Delta E^H=\Delta E$ for the blocks in~\eqref{eq:M_N_def}.

The corresponding block-structure- and symmetry-structure-preserving eigenpair backward errors
$\eta^{\mathcal B}(J,E,\lambda,x)$ and $\eta^{\mathcal S}(J,E,\lambda,x)$ are defined by~\eqref{eq:def_JE_eperr}
and~\eqref{eq:def_strJE_eperr}, respectively. We first discuss under which conditions these backward errors
are finite.

\begin{proposition}\label{prop1:14.9.17}\rm
Let $L(z)$ be a pencil as in~\eqref{eq:defL_cz}, and let $\lambda \in\mathbb C$ and $x=[x_1^T~x_2^T~x_3^T]^T$
be such that $x_1,\,x_2 \in \C^{n}$ and $x_3\in \C^m$.
Then for any $\Delta_J,\,\Delta_E \in \C^{n , n}$ and corresponding
$\Delta L(z)=\Delta_M +z \Delta_N \in \mathcal B$, we have
$(L-\Delta L)(\lambda)x=0$ if and only if
\begin{eqnarray}
(\Delta_J + \lambda \Delta_E)x_2 &=&(J-R+\lambda E)x_2 +B x_3, \label{equi_1}\\
(\Delta_J+\lambda \Delta_E)^H x_1&=& (-J-R-\lambda E)x_1,\label{equi_2}\\
0&=& B^Hx_1+Sx_3, \label{equi_3}
\end{eqnarray}
i.~e., $\eta^{\mathcal B}(J,E,\lambda,x)$ is finite if and only if there exists matrices
$\Delta_J$ and $\Delta_E$ such that these equations are satisfied.
\end{proposition}

In the next lemma we
present conditions that are equivalent to the existence of matrices $\Delta_J$ and
$\Delta_E$ that satisfy the first two of the three equations in Proposition~\ref{prop1:14.9.17}.

\begin{lemma}\label{thm;equivalentcond}
Let $L(z)$ be a pencil as in~\eqref{eq:defL_cz}, and let $\lambda \in i\R$ and $x=[x_1^T~x_2^T~x_3^T]^T$
be such that $x_1,\,x_2 \in \C^{n}$ and $x_3\in \C^m$. Furthermore, set
$r:=(J-R+\lambda E)x_2+Bx_3$ and $s:=(-J-R-\lambda E)x_1$. Then the following statements are equivalent.
\begin{enumerate}
 \item There exist $\Delta_J \in \C^{n , n}$ and $\Delta_E \in \C^{n , n}$
 satisfying~\eqref{equi_1} and~\eqref{equi_2}.
 \item There exist $ \Delta \in \C^{n , n}$ such that $ \Delta x_2=r$ and $ \Delta^H x_1=s$.
 \item The identity $x_3^HB^Hx_1=0$ is satisfied.
\end{enumerate}
Moreover, we have
\begin{eqnarray}\label{JEequi_4}
\inf\left\{ {\|\Delta_J\|}_F^2+{\|\Delta_E\|}_F^2 ~\Big|~ \Delta_J,\,\Delta_E \in \C^{n , n}
~{\rm satisfy}~\eqref{equi_1}\,{\rm and}~\eqref{equi_2}  \right \}\qquad\qquad \nonumber \\
=\inf \left\{\left.\frac{{\| \Delta\|}_F^2}{1+|\lambda|^2}~\right|~  \Delta \in \C^{n , n},\,
\Delta x_2=r,~ \Delta ^H x_1=s\right \}.
\end{eqnarray}
\end{lemma}

\proof ``$1)\Rightarrow 2)$'': Let $\Delta_J \in \C^{n , n}$ and $\Delta_E \in \C^{n , n}$ be such
that~\eqref{equi_1} and~\eqref{equi_2} are satisfied. Then by setting $\Delta = \Delta_J + \lambda \Delta_E$
we get $\Delta x_2=r$, $\Delta^H x_1=s$ which shows 2). Furthermore, using
the Cauchy-Schwarz inequality (in $\mathbb R^2$), we obtain
\[
{\|\Delta\|}_F^2 \leq \big(\|\Delta_J \|_F + |\lambda|\,\|\Delta_E\|_F\big)^2
\leq (1+|\lambda|^2) ({\|\Delta_J\|}_F^2+{\|\Delta_E\|}_F^2).
\]
This implies
\begin{eqnarray*}
  \inf \left\{\left.\frac{{\|\Delta\|}_F^2}{1+|\lambda|^2}~\right|~
  \Delta_J,\,\Delta_E \in \C^{n , n} ~{\rm satisfy}~\eqref{equi_1}\,{\rm and}~\eqref{equi_2},~
  \Delta=\Delta_J+\lambda \Delta_E\right \} \qquad\qquad\\
\leq \inf\left\{ {\|\Delta_J\|}_F^2+{\|\Delta_E\|}_F^2 ~\big|~ \Delta_J,\,\Delta_E \in \C^{n , n} ~{\rm satisfy}~
  \eqref{equi_1}\,{\rm and}~\eqref{equi_2}  \right \},
\end{eqnarray*}
and thus
\begin{eqnarray}\label{equi_5}
&&  \inf \left\{\left.\frac{{\|\Delta\|}_F^2}{1+|\lambda|^2}~\right|~ \Delta \in \C^{n , n},\,
  \Delta x_2=r,~\Delta ^H x_1=s\right \}\nonumber \\
&&\qquad\qquad\leq\inf\left\{ {\|\Delta_J\|}_F^2+{\|\Delta_E\|}_F^2 ~\big|~ \Delta_J,\,
\Delta_E \in \C^{n , n} ~{\rm satisfy}~\eqref{equi_1}\,{\rm and}~\eqref{equi_2}  \right \}\qquad\qquad
\end{eqnarray}
which yields ``$\geq$'' in~\eqref{JEequi_4}.

\noindent
``$2) \Rightarrow 1)$'': Conversely, suppose that $\Delta \in \C^{n , n} $ satisfies $\Delta x_2=r$
and $\Delta ^H x_1=s$. Then by setting $\Delta_J=\frac{\Delta}{1+|\lambda|^2}$ and
$\Delta_E=\frac{\bar \lambda \Delta}{1+|\lambda|^2}$ we get $\Delta_J+\lambda \Delta_E =\Delta$
and hence, $\Delta_J$ and $\Delta_E$ satisfy~\eqref{equi_1} and~\eqref{equi_2} which proves 1).
Furthermore, we obtain
\[
{\|\Delta_J\|}_F^2+{\|\Delta_E\|}_F^2 = \frac{{\|\Delta \|}_F^2}{(1+|\lambda|^2)^2}
+\frac{|\lambda|^2{\|\Delta \|}_F^2}{(1+|\lambda|^2)^2}=\frac{{\|\Delta \|}_F^2}{1+|\lambda|^2}.
\]
This implies
\begin{eqnarray*}
  \inf\left\{{\|\Delta_J\|}_F^2+{\|\Delta_E\|}_F^2 \,\left|\,  \Delta \in \C^{n , n},\, \Delta x_2=r,\,\Delta ^H x_1=s,\,
 \Delta_J=\frac{\Delta}{1+|\lambda|^2},\Delta_E =\frac{\bar \lambda \Delta}{1+|\lambda|^2}\right.\right\}\qquad \\
  =\inf \left\{\frac{{\|\Delta\|}_F^2}{1+|\lambda|^2}~\Big |~ \Delta \in \C^{n , n},\,
  \Delta x_2=r,~\Delta ^H x_1=s\right \}
\end{eqnarray*}
and hence
\begin{eqnarray*}
  \inf\left\{ \left.{\|\Delta_J\|}_F^2+{\|\Delta_E\|}_F^2 ~\right|~ \Delta_J,\,\Delta_E \in \C^{n , n}
  ~{\rm satisfy}~\eqref{equi_1}\,{\rm and}~\eqref{equi_2}  \right \} \qquad \\
  \leq\inf \left\{\frac{{\|\Delta\|}_F^2}{1+|\lambda|^2}~\Big |~ \Delta \in \C^{n , n},\,
  \Delta x_2=r,~\Delta ^H x_1=s\right \}
\end{eqnarray*}
which proves ``$\leq$'' in~\eqref{JEequi_4}.

\noindent
``$2)\Leftrightarrow 3)$'': This follows from Theorem~\ref{thm:spcl_matrix_map}, because there exists
 $\Delta \in \C^{n , n}$ satisfying $\Delta x_2=r$ and $\Delta^H x_1=s $
 if and only if $x_2^Hs=r^Hx_1$. Since $\lambda$ is purely imaginary, this latter equation
 is equivalent to $x_3^HB^Hx_1=0$.
\eproof

\begin{theorem}\label{thm:ustrJE}
Let $L(z)$ be a pencil as in~\eqref{eq:defL_cz}, and let $\lambda \in i \R$ and $x \in \C^{2n+m}\setminus \{0\}$.
Partition $x=[x_1^T~x_2^T~x_3^T]^T$ such that $x_1,\,x_2\in \C^n $ and $x_3 \in \C^{m}$ and set
$r=(J-R+\lambda E)x_2+Bx_3$ and $s=-(J+R+\lambda E)x_1$.
Then $\eta^{\mathcal B}(J,E,\lambda,x)$ is finite if and only if $x_3=0$ and $B^Hx_1=0$.
If the latter conditions hold then
\begin{eqnarray}\label{ustrbacerr_JE}
 \eta^{\mathcal B}(J,E,\lambda,x)~ = ~ \frac{{\|\widehat {\Delta}\|}_F}{\sqrt{1+|\lambda|^2}}\quad
 \text{and}\quad \eta^{\mathcal B}(J,E,\lambda)=\frac{\sigma_{\min}(J-R+\lambda E)}{\sqrt{1+|\lambda|^2}},
\end{eqnarray}
where $\widehat{\Delta}$ is given by
\begin{equation*}
\widehat{\Delta}= \left\{\begin{array}{ll}
\frac{rx_2^H}{\|x_2\|^2} & \mbox{ if } x_1 = 0, \\[1ex]
\frac{x_1s^H}{\|x_1\|^2} & \mbox{ if } x_2 = 0,\\[1ex]
\frac{rx_2^H}{\|x_2\|^2} + \frac{x_1s^H}{\|x_1\|^2}\left(I_{n} - \frac{x_2x_2^H}{\|x_2\|^2}\right)
& \mbox{ otherwise. }  \end{array}\right.
\end{equation*}
\end{theorem}

\proof
Combining Proposition~\ref{prop1:14.9.17} and Lemma~\ref{thm;equivalentcond},
we obtain that $\eta^{\mathcal B}(J,E,\lambda,x)$ is finite if and only if $x$ satisfies $x_3^HB^Hx_1=0$ and
$B^Hx_1+Sx_3=0$, or equivalently,  $x_3=0$ and $B^Hx_1=0$, since $S$ is definite.
Thus, assume that $x$ satisfies $x_3=0$ and $B^Hx_1=0$. Then we obtain
\begin{eqnarray}\label{eq:JE23aug1}\nonumber
\eta^{\mathcal B}(J,E,\lambda,x)&=&\inf\Big\{\|[\Delta_J~\Delta_E]\|_{F} \;\Big |
\;\big((M-\Delta_M)+\lambda(N-\Delta_N)\big)x=0,\;\Delta_M + z \Delta N \in \mathcal B\Big\}\\
&=&\inf\Big\{\|[\Delta_J~\Delta_E]\|_{F}\, \Big |
\Delta_J,\,\Delta_E \in \C^{n , n}~{\rm satisfy}~\eqref{equi_1}~{\rm and}~\eqref{equi_2} \Big\}\nonumber\\
&=&\inf \left\{\left.\frac{{\|\Delta\|}_F}{\sqrt{1+|\lambda|^2}}\,\right|~
  \Delta \in \C^{n , n},\, \Delta x_2=r,~\Delta ^H x_1=s\right \}\nonumber\\
&=&\frac{1}{\sqrt{1+|\lambda|^2}}\inf \Big\{{\|\Delta\|}_F\, \Big|~
  \Delta \in \C^{n , n},\, \Delta x_2=r,~\Delta ^H x_1=s\Big\},
\end{eqnarray}
where the second last equality is due to Lemma~\ref{thm;equivalentcond}.
Thus the formula for $\eta^{\mathcal B}(J,E,\lambda,x)$ in~\eqref{ustrbacerr_JE} follows from
Theorem~\ref{thm:spcl_matrix_map} for the case $x_1,x_2\neq 0$ and for the case $x_1=0$ or $x_2=0$ it
is straightforward. (Indeed, in the case $x_1=0$ any matrix $\Delta$ with $\Delta x_2=r$ satisfies
$\|\Delta\|_F\geq\frac{\|r\|}{\|x_2\|}$ and $\widehat\Delta=\frac{rx_2^H}{\|x_2\|^2}$ is a matrix
for which equality is attained. The case $x_2=0$ is analogous.)

Next we will prove the formula for $\eta^{\mathcal B}(J,E,\lambda)$ in~\eqref{ustrbacerr_JE}. To this end, let
$$\mathcal M:=\big\{y=\mat{ccc}y_1^T&y_2^T&0\rix\in\mathbb C^{2n+m}\,\big|\,y_1,y_2\in\mathbb C^{n},\,
(y_1,y_2)\neq(0,0),\,B^Hy_1=0\big\}.$$ Then we obtain
\begin{eqnarray}\label{eq:JE23aug2}
&&\big(\sqrt{1+|\lambda|^2}\big)\cdot\eta^{\mathcal B}(J,E,\lambda) =
\big(\sqrt{1+|\lambda|^2}\big) \cdot\inf_{y \in \C^{2n+m}\setminus \{0\}}\eta^{\mathcal B}(J,E,\lambda,y)\nonumber\\
&=& \inf_{y\in\mathcal M}\inf \Big\{{\|\Delta\|}_F\, \Big|~
\Delta \in \C^{n , n},\, \Delta y_2=(J-R+\lambda E)y_2,~ \Delta ^H y_1=-(J+R+\lambda E)y_1 \Big \} \nonumber \\
&\geq& \inf_{y\in\mathcal M}\inf \Big\{{\|\Delta\|}\, \Big|~  \Delta \in \C^{n , n},\,
\Delta y_2=(J-R+\lambda E)y_2,~ \Delta ^H y_1=-(J+R+\lambda E)y_1\Big \},\qquad
\end{eqnarray}
where the second equality is due to~\eqref{eq:JE23aug1} and the inequality in the last line follows from the fact
that for any $\Delta \in \C^{n,n}$, we have $\|\Delta\|\leq {\|\Delta\|}_F$. Defining
\begin{eqnarray*}
\mu := \inf_{y\in\mathcal M}
\inf \Big\{{\|\Delta\|}\, \Big|~  \Delta \in \C^{n , n},\,
\Delta y_2=(J-R+\lambda E)y_2,~
 \Delta ^H y_1=-(J+R+\lambda E)y_1\Big \},
\end{eqnarray*}
we get by applying Theorem~\ref{thm:spcl_matrix_map} for the case of the spectral norm that
\begin{eqnarray}\label{eq:JE23aug33}
\mu &=& \inf_{y\in\mathcal M}
\max\left\{\frac{\|(J-R+\lambda E)^H y_1\|}{\|y_1\|},\frac{\|(J-R+\lambda E)y_2\|}{\|y_2\|}\right\}\\
&=&\min\left\{\inf_{y_1 \in \C^{n}\setminus \{0\},B^Hy_1=0}\frac{\|(J-R+\lambda E)^H y_1\|}{\|y_1\|},
\inf_{y_2 \in \C^{n}\setminus \{0\}}\frac{\|(J-R+\lambda E)y_2\|}{\|y_2\|}\right\},\qquad\label{eq:JE23aug3}
\end{eqnarray}
where in~\eqref{eq:JE23aug33} we interpret the undefined expressions $\frac{0}{0}$ that occur in the
cases $y_1=0$ or $y_2=0$ as being equal to zero.
Let the columns of $U=[u_1,~\ldots,~u_k]\in\mathbb C^{n,k}$ form
an orthonormal basis of $\operatorname{null}(B^H)$. Then
\begin{eqnarray}\label{eq:JE23aug4}
\inf_{y_1\in \C^{n}\setminus\{0\},B^Hy_1=0}\frac{{\|(J-R+\lambda E)^Hy_1\|}^2}{{\|y_1\|}^2}
&=&\inf_{y_1\in\operatorname{null}(B^H)\setminus\{0\}}\frac{{\|(J-R+\lambda E)^Hy_1\|}^2}{{\|y_1\|}^2} \nonumber \\
&=& \inf_{\alpha \in \C^{k}\setminus\{0\}} \frac{{\|(J-R+\lambda E)^HU \alpha\|}^2}{{\|\alpha\|}^2} \nonumber\\
&=& \Big(\sigma_{\min}\left((J-R+\lambda E)^HU\right)\Big)^2.
\end{eqnarray}
By inserting~\eqref{eq:JE23aug4} in~\eqref{eq:JE23aug3}, we get
\begin{eqnarray}\label{eq:JE23aug5}
\mu &=& \min\left\{\sigma_{\min}\left((J-R+\lambda E)^H\right),\,\sigma_{\min}\left((J-R+\lambda E)^HU\right)\right\}\nonumber \\
&=& \sigma_{\min}\left((J-R+\lambda E)^H\right)=\sigma_{\min}\big((J-R+\lambda E)\big),
\end{eqnarray}
Using the value of $\mu$ from~\eqref{eq:JE23aug5}, we show that equality holds in~\eqref{eq:JE23aug2} by constructing
$\Delta$ such that $\|\Delta\|={\|\Delta\|}_F=\mu$. For this, let $u$ and $v$ respectively be unit left and  right singular vectors of $(J-R+\lambda E)$ corresponding to the singular value
$\sigma^*:=\sigma_{\min}\left((J-R+\lambda E)\right)$ and consider $\widetilde \Delta:=\sigma^*uv^H$.
Then, clearly $\|\widetilde \Delta\|={\|\widetilde \Delta\|}_F=\sigma^*$ as $\widetilde \Delta$ is of rank one, and
\[
\widetilde \Delta v=\sigma^*u =(J-R+\lambda E)v\quad \text{and}\quad \widetilde \Delta^Hu=\sigma^*v=(J-R+\lambda E)^H u.
\]
Thus we have equality in~\eqref{eq:JE23aug2}, i.e.,
\[
\eta^{\mathcal B}(J,E,\lambda)=\frac{\mu}{\sqrt{1+|\lambda|^2}}=\frac{{\|\widetilde \Delta\|}_F}{\sqrt{1+|\lambda|^2}}
=\frac{\sigma_{\min}(J-R+\lambda E)}{\sqrt{1+|\lambda|^2}}
\]
which finishes the proof. \eproof
\medskip
Next we aim to compute the symmetry-structure-preserving eigenpair error $\eta^{\mathcal S}(J,E,\lambda,x)$, i.e.,
when we have $\Delta_J^H=-\Delta_J$ and $\Delta_E^H=\Delta_E$ in the pencils $L(z)=\Delta_M+z\Delta_N\in\mathcal S$.
We start with a criterion for the finiteness of the eigenpair error, where we focus on the case that
$\lambda$ is on the imaginary axis.

\begin{remark}\label{rem11:15.9.17}\rm
Let $L(z)$ be a pencil as in~\eqref{eq:defL_cz}, and let $\lambda \in i\mathbb R$ and $x=[x_1^T~x_2^T~x_3^T]$
be such that $x_1,\,x_2 \in \C^{n}$ and $x_3\in \C^m$.
Then using $\Delta_J^H=-\Delta_J$ and $\Delta_E^H=\Delta_E$ and also the fact that $\lambda$
is purely imaginary, the equations \eqref{equi_1}--\eqref{equi_3} take the form
\begin{eqnarray*}
(\Delta_J + \lambda \Delta_E)x_2 &=&(J-R+\lambda E)x_2 +B x_3, \\
(-\Delta_J-\lambda \Delta_E) x_1 &=& (-J-R-\lambda E)x_1,\\
0&=& B^Hx_1+Sx_3.
\end{eqnarray*}
Thus, combining the first two of these equations, we find that $\eta^{\mathcal S}(J,E,\lambda,x)$
is finite if and only if there exist $\Delta_J\in\text{SHerm}(n)$ and $\Delta_E\in\text{Herm}(n)$
such that the equations
\begin{eqnarray}\label{eq:equi_strJE1}
 (\Delta_J+\lambda \Delta_E)\mat{cc}x_2&x_1\rix&=&\mat{cc}(J-R+\lambda E)x_2 +B x_3&
 (J+R+\lambda E)x_1\rix,\\
 0&=& B^Hx_1+Sx_3 \nonumber
\end{eqnarray}
are satisfied.
\end{remark}\\
We start with a lemma that contains equivalent conditions for equation~\eqref{eq:equi_strJE1}
to be satisfied.
\begin{lemma} \label{lem1:str_JE}
 Let $L(z)$ be a pencil as in~\eqref{eq:defL_cz}, and let
$\lambda \in i \R$ and $x \in \C^{2n+m}\setminus \{0\}$. Partition $x=[x_1^T~x_2^T~x_3^T]^T$
such that $x_1,\,x_2 \in \C^{n}$ and $x_3 \in \C^{m}$, and define
\[
X=\mat{cc}x_2&x_1\rix \quad \text{and}\quad Y=\mat{cc}(J-R+\lambda E)x_2 +B x_3&(J+R+\lambda E)x_1\rix.
\]
Then the following statements are equivalent.
\begin{enumerate}
 \item There exist $\Delta_J \in {\rm SHerm}(n)$ and $\Delta_E \in {\rm Herm}(n)$
 satisfying~\eqref{eq:equi_strJE1}.
 \item There exist $\Delta \in {\rm SHerm}(n)$ such that $\Delta X=Y$.
 \item $X$ and $Y$ satisfy $Y^HX=-X^HY$ and $YX^{\dagger}X=Y$.
\end{enumerate}
Moreover, we have
\begin{eqnarray}\label{eq:JE_normequality}
  \inf\left\{\left. {\|\Delta_J\|}_F^2+{\|\Delta_E\|}_F^2 ~\right|~ \Delta_J \in {\rm SHerm}(n),\,\Delta_E \in
  {\rm Herm}(n)~{\rm satisfying}~\eqref{eq:equi_strJE1}  \right \} \qquad\nonumber \\
  =\inf \left\{\left.\frac{{\|\Delta\|}_F^2}{1+|\lambda|^2}~\right|~ \Delta \in {\rm SHerm}(n),\,
  \Delta X=Y\right \}.
 \end{eqnarray}
\end{lemma}

\proof ``$1)\Rightarrow 2)$'': Let $\Delta_J \in {\rm SHerm}(n)$ and $\Delta_E \in {\rm Herm}(n)$
be such that they satisfy~\eqref{eq:equi_strJE1}, then by setting $\Delta = \Delta_J + \lambda \Delta_E$
we get $\Delta X=Y$, and $\Delta \in {\rm SHerm}(n) $ as $\lambda \in i\R$. The inequality
``$\geq$'' in~\eqref{eq:JE_normequality} then follows by the same  arguments as in ``$1)\Rightarrow 2)$'' in
the proof of Lemma~\ref{thm;equivalentcond}.

\noindent
``$2)\Rightarrow 1)$'': Conversely, let $\Delta \in {\rm SHerm}(n)$ be such that $\Delta X=Y$. Then setting
\[
\Delta_J=\frac{\Delta}{1+|\lambda|^2}\quad\mbox{and}\quad\Delta_E=\frac{\bar \lambda \Delta}{1+|\lambda|^2},
\]
we obtain $(\Delta_J+\lambda \Delta_E)X=Y$ as well as
$\Delta_J \in {\rm SHerm}(n)$ and $\Delta_E \in {\rm Herm}(n)$, since $\lambda \in i\R$.
Again, the proof ``$\leq$'' in~\eqref{eq:JE_normequality} follows by arguments similar to those
in the part ``$2)\Rightarrow 1) $'' in the proof of Lemma~\ref{thm;equivalentcond}.

\noindent
``$2)\Leftrightarrow 3)$'': This follows immediately by Theorem~\ref{eq:SHerm_map}.
\eproof

\begin{theorem}\label{thm:strJE}
Let $L(z)$ be a pencil defined as in~\eqref{eq:defL_cz}, and let
$\lambda \in i \R$ and $x \in \C^{2n+m}\setminus \{0\}$. Partition $x=[x_1^T~x_2^T~x_3^T]^T$
 such that $x_1,\,x_2 \in \C^{n}$ and $x_3 \in \C^{m}$ and set
 \[
X=\mat{cc}x_2&x_1\rix \quad \text{and}\quad Y=\mat{cc}(J-R+\lambda E)x_2 +B x_3&(J+R+\lambda E)x_1\rix.
\]
Then $\eta^{\mathcal S}(J,E,\lambda,x)$ is finite if and only if
 $Y^HX=-X^HY$, $YX^{\dag}X=Y$ and $B^Hx_1+Sx_3=0$. If the three latter conditions are satisfied, then
\begin{equation}\label{eq:strbacerr_JE}
  \eta^{{\mathcal S}}(J,E,\lambda,x)=\sqrt{\frac{1}{1+|\lambda|^2}
  \left(2{\|YX^{\dag}\|}_F^2-\operatorname{trace}\big(YX^{\dag}(YX^\dag)^HXX^\dag\big)\right)}.
\end{equation}
\end{theorem}

\proof
Combining Remark~\ref{rem11:15.9.17} and Lemma~\ref{lem1:str_JE} it follows that $\eta^{\mathcal S}(J,E,\lambda,x)$
is finite if and only if $x$ satisfies
\[
Y^HX=-X^HY,\quad YX^{\dag}X=Y\quad {\rm and}\quad B^Hx_1+Sx_3=0.
\]
In the following let us assume that these conditions on $x$ are satisfied. Then we obtain
\begin{eqnarray*}
 \eta^{\mathcal S}(J,E,\lambda,x)&=& \inf\Big\{\|[\Delta_J~\Delta_E]\|_{F} \Big |\,
 \Delta_J \in {\rm SHerm}(n),\,\Delta_E \in {\rm Herm}(n)~{\rm satisfy}~\eqref{eq:equi_strJE1} \Big\}\\
&=& \frac{1}{\sqrt{1+|\lambda|^2}}\cdot\inf \Big\{{\|\Delta\|}_F~\Big |~ \Delta \in {\rm SHerm}(n),\,
  \Delta X=Y\Big \},
\end{eqnarray*}
where the last equality is due to Lemma~\ref{lem1:str_JE}. Hence~\eqref{eq:strbacerr_JE} follows by using
Theorem~\ref{eq:SHerm_map}.
\eproof

\subsection{Perturbations only in $R$ and $E$}
In this section, we consider the case that only the blocks $R$ and $E$ in a pencil $L(z)$ as
in~\eqref{eq:defL_cz} are perturbed. Let $\lambda\in \C$ and $x \in \C^{2n+m}\setminus \{0\}$.
Then by the terminology outlined in Section~\ref{sec:section_def_epbr}, the block- and
symmetry-structure-preserving eigenpair backward errors
$\eta^{\mathcal B}(R,E,\lambda,x)$ and $\eta^{\mathcal S}(R,E,\lambda,x)$ are defined by
\begin{equation}\label{eq:def_RE_eperr}
\eta^{\mathcal B}(R,E,\lambda,x)=\inf\Big\{\|[\Delta_R~\Delta_E]\|_{F} \, \Big |\,
\big((M-\Delta_M)+\lambda(N-\Delta_N)\big)x=0,\,\Delta_M + z \Delta_N \in \mathcal B\Big\},
\end{equation}
and
\begin{equation}\label{eq:deff_strRE_eperror}
\eta^{\mathcal S}(R,E,\lambda,x)=\inf\Big\{\|[\Delta_R~\Delta_E]\|_{F} \, \Big |\,
\big((M-\Delta_M)+\lambda(N-\Delta_N)\big)x=0,\,\Delta_M + z \Delta_N \in \mathcal S\Big\},
\end{equation}
respectively, where $\mathcal B$ is the set of all pencils
of the form $\Delta L(z)=\Delta_M + z \Delta_N$ with the block-structure
\begin{equation}\label{14.9.17:pencil}
\Delta_M= \mat{ccc} 0 & -\Delta_R & 0\\ -\Delta _R^H &0&0 \\0&0&0\rix \quad \text{and} \quad
\Delta_N= \mat{ccc} 0 & \Delta_E & 0\\ -\Delta_E^H &0&0 \\0&0&0\rix.
\end{equation}
and $\Delta_R,\Delta_E\in\mathbb C^{n,n}$, while $\mathcal S$ is the corresponding set
of pencils $\Delta L(z)=\Delta_M + z \Delta_N$ as in~\eqref{14.9.17:pencil} with
$\Delta_R,\Delta_E\in\text{Herm}(n)$.

We highlight that in the case that only the block-structure is preserved, the perturbation
matrices in~\eqref{14.9.17:pencil} have exactly the same structure as the ones in~\eqref{eq:M_N_def}
and hence by following exactly the same lines as in the previous section, we obtain the following
theorem which shows that the values
of $\eta^{\mathcal B}(R,E,\lambda,x)$, and also of $\eta^{\mathcal B}(R,E,\lambda):=
\inf_{x \in \C^{2n+m}\setminus \{0\}}\eta^{\mathcal B}(R,E,\lambda,x)$ are equal to the corresponding
values $\eta^{\mathcal B}(J,E,\lambda,x)$ and $\eta^{\mathcal B}(J,E,\lambda)$.

\begin{theorem}\label{thm:ustrRE}
Let $L(z)$ be a pencil as in~\eqref{eq:defL_cz}, and let
$\lambda \in i \R$ and $x \in \C^{2n+m}\setminus \{0\}$. Partition $x=[x_1^T~x_2^T~x_3^T]^T$
such that $x_1,\,x_2\in \C^n$ and $x_3 \in \C^{m}$ and set
$r=(J-R+\lambda E)x_2 +Bx_3$ and $s=-(J+R+\lambda E)x_1.$
Then $\eta^{\mathcal B}(R,E,\lambda,x)$ is finite if and only if $x_3=0$ and $B^Hx_1=0$.
If the latter conditions hold then
\begin{eqnarray*}\label{eq:RE_unsstrberr}
 \eta^{\mathcal B}(R,E,\lambda,x)=\eta^{\mathcal B}(J,E,\lambda,x)~ = ~
\frac{{\|\widehat {\Delta}\|}_F}{\sqrt{1+|\lambda|^2}},
\end{eqnarray*}
and
\begin{eqnarray*}
\eta^{\mathcal B}(R,E,\lambda)=\eta^{\mathcal B}(J,E,\lambda)= \frac{\sigma_{\min}(J-R+\lambda E)}{\sqrt{1+|\lambda|^2}},
\end{eqnarray*}
where $\widehat{\Delta}$ is given by
\begin{equation*}
\widehat{\Delta}= \left\{\begin{array}{ll}
\frac{rx_2^H}{\|x_2\|^2} & \mbox{ if } x_1 = 0, \\[1ex]
\frac{x_1s^H}{\|x_1\|^2} & \mbox{ if } x_2 = 0,\\[1ex]
\frac{rx_2^H}{\|x_2\|^2} + \frac{x_1s^H}{\|x_1\|^2}\left(I_{n} - \frac{x_2x_2^H}{\|x_2\|^2}\right)
& \mbox{ otherwise. }  \end{array}\right.
\end{equation*}
\end{theorem}

\proof The proof follows exactly the same as the proof of Theorem~\ref{thm:ustrJE} by
just replacing $\Delta J$ with $-\Delta R$.
\eproof

Next, we turn to the eigenpair backward error $\eta^{{\mathcal S}}(R,E,\lambda,x)$ for purely imaginary
$\lambda \in i \R$ and
$x=[x_1^T~x_2^T~x_3^T]^T\in \C^{2n+m}\setminus \{0\}$. Note that in this case $\Delta_R^H=\Delta_R$
and $\Delta_E^H=\Delta_E$. In particular, the perturbations now have a different symmetry structure
than the corresponding ones from the previous section, so that we expect the backward error
$\eta^{{\mathcal S}}(R,E,\lambda,x)$ to differ from $\eta^{{\mathcal S}}(J,E,\lambda,x)$.
We start again with a criterion for the finiteness of $\eta^{{\mathcal S}}(R,E,\lambda,x)$ and
continue with a lemma giving equivalent conditions.

\begin{remark}\label{rem1:15.9.17}\rm
Let $L(z)$ be a pencil as in~\eqref{eq:defL_cz}, and let $\lambda \in\mathbb C$ and $x=[x_1^T~x_2^T~x_3^T]$
be such that $x_1,\,x_2 \in \C^{n}$ and $x_3\in \C^m$.
Then using $\Delta_R^H=-\Delta_R$ and $\Delta_E^H=\Delta_E$ and also the fact that $\lambda$
is purely imaginary, we find that there exist $\Delta_R,\,\Delta_E \in \text{Herm(n)}$ and correspondingly
$\Delta L(z)=\Delta_M +z \Delta_N \in \mathcal S$ such that $(L-\Delta L)(\lambda)x=0$ if and only if
\begin{eqnarray}
(-\Delta_R + \lambda \Delta_E)x_2 &=&(J-R+\lambda E)x_2 +B x_3 \label{REequi_1}\\
(-\Delta_R+\lambda \Delta_E)^H x_1 &=& (-J-R-\lambda E)x_1\label{REequi_2}\\
0&=& B^Hx_1+Sx_3 \label{REequi_3}.
\end{eqnarray}
Thus, $\eta^{{\mathcal S}}(R,E,\lambda,x)$ is finite if and only if there exist $\Delta_R,\,\Delta_E \in \text{Herm(n)}$
such that~\eqref{REequi_1}--\eqref{REequi_3} are satisfied.
\end{remark}

\begin{lemma}\label{lem1:str_RE}
Let $L(z)$ be a pencil as in~\eqref{eq:defL_cz}, and let
$\lambda \in i \R$ and $x \in \C^{2n+m}\setminus \{0\}$. Partition $x=[x_1^T~x_2^T~x_3^T]^T$
 such that $x_1,\,x_2 \in \C^n$ and $x_3 \in \C^{m}$ and let $r=(J-R+\lambda E)x_2 +B x_3$
and $s=(-J-R-\lambda E)x_1$. Then the following statements are equivalent.
\begin{enumerate}
\item  There exist $\Delta_R,\, \Delta_E \in {\rm Herm}(n)$ satisfying~\eqref{REequi_1} and~\eqref{REequi_2}.
\item  There exist $\Delta \in \C^{n , n}$ such that $\Delta x_2=r$ and $\Delta ^H x_1=s$.
\item  The identity  $x_3^HB^Hx_1=0$ is satisfied.
\end{enumerate}
Moreover, we have
\begin{eqnarray}\label{REequi_4}
&&\inf\left\{ \left.{\|\Delta_R\|}_F^2+{\|\Delta_E\|}_F^2 ~\right|~ \Delta_R,\,\Delta_E \in
{\rm Herm}(n)~{\rm satisfy}~\eqref{REequi_1}\,{\rm and}~\eqref{REequi_2}  \right \}\nonumber\\
&&\quad =
\inf \left\{\left.{\left\|\frac{\Delta +\Delta^H}{2}\right\|}_F^2+
\frac{1}{|\lambda|^2}{\left\|\frac{\Delta -\Delta^H}{2}\right\|}_F^2~\right|~ \Delta \in \C^{n , n},\,
  \Delta x_2=r,~\Delta^H x_1=s\right \}.\qquad
\end{eqnarray}
\end{lemma}

\proof
``$1) \Rightarrow 2)$'': Let $\Delta_R,\,\Delta_E \in {\rm Herm}(n)$ be such that
 they satisfy~\eqref{REequi_1} and~\eqref{REequi_2}. Setting $\Delta = -\Delta_R + \lambda \Delta_E$, we get
 $\Delta x_2=r,~\Delta^Hx_1=s$. Also note that $-\Delta_R$ and
 $\lambda \Delta_E$ are the unique Hermitian and skew-Hermitian
parts of $\Delta$, respectively, i.e.,  $\Delta_R=-(\Delta +\Delta^H)/2$ and
$\lambda \Delta_E=(\Delta-\Delta^H)/2 $.
This implies
\[
{\|\Delta_R\|}_F^2+{\|\Delta_E\|}_F^2={\left\|\frac{\Delta +\Delta^H}{2}\right\|}_F^2+
\frac{1}{|\lambda|^2}{\left\|\frac{\Delta -\Delta^H}{2}\right\|}_F^2
\]
 and
\begin{eqnarray*}
&&\inf\left\{\left. {\|\Delta_R\|}_F^2+{\|\Delta_E\|}_F^2 ~\right|~ \Delta_R,\,\Delta_E \in {\rm Herm}(n) ~{\rm satisfy}~\eqref{REequi_1}\,{\rm and}~\eqref{REequi_2} \right \}
\\
&&\qquad=\inf \Bigg\{{\left\|\frac{\Delta +\Delta^H}{2}\right\|}_F^2+
\frac{1}{|\lambda|^2}{\left\|~\frac{\Delta -\Delta^H}{2}\right\|}_F^2\,\Bigg |~
\Delta=-\Delta_R+\lambda \Delta_E,~ \Delta_R,\,\Delta_E \in {\rm Herm}(n)\qquad \\
&&\phantom{\qquad=\inf \Bigg\{{\left\|\frac{\Delta +\Delta^H}{2}\right\|}_F^2+
\frac{1}{|\lambda|^2}{\left\|~\frac{\Delta -\Delta^H}{2}\right\|}_F^2\,\Bigg |~\Delta=-\Delta_R+\lambda}
~{\rm satisfy}~\eqref{REequi_1}\,{\rm and}~\eqref{REequi_2}
\Bigg \}.
\end{eqnarray*}
Thus, we obtain
\begin{eqnarray}\label{REequi_5}
&& \inf\left\{ {\|\Delta_R\|}_F^2+{\|\Delta_E\|}_F^2 ~\big|~ \Delta_J,\,\Delta_E \in
 {\rm Herm}(n)~{\rm satisfy}~\eqref{REequi_1}\,{\rm and}\,\eqref{REequi_2}  \right \}\\
&&\quad\geq
  \inf \left\{\left.{\left\|\frac{\Delta +\Delta^H}{2}\right\|}_F^2+
\frac{1}{|\lambda|^2}{\left\|\frac{\Delta -\Delta^H}{2}\right\|}_F^2~\right|~ \Delta \in \C^{n , n},\,
 \Delta x_2=r,~\Delta ^H x_1=s\right \} \qquad \nonumber
  \end{eqnarray}
which gives ``$\geq$'' in~\eqref{REequi_4}.

\noindent
 ``$2) \Rightarrow 1)$'': Suppose that $\Delta \in \C^{n , n} $ is such that $\Delta x_2=r$
 and $\Delta ^H x_1=s$. Then, by setting $\Delta_R=-\frac{\Delta+\Delta^H}{2}$ and
 $\Delta_E=\frac{\bar \lambda}{|\lambda|^2}(\frac{ \Delta-\Delta^H}{2})$, we get
 $\Delta_R,\,\Delta_E\in {\rm Herm}(n)$ such that~\eqref{REequi_1} and~\eqref{REequi_2}
 are satisfied, because $-\Delta_R+\lambda \Delta_E =\Delta$. Also, we have
\[
{\|\Delta_R\|}_F^2+{\|\Delta_E\|}_F^2 = {\left\|\frac{\Delta +\Delta^H}{2}\right\|}_F^2+
\frac{1}{|\lambda|^2}{\left\|\frac{\Delta -\Delta^H}{2}\right\|}_F^2
\]
which implies
\begin{eqnarray*}
&&\inf \left\{\left.{\left\|\frac{\Delta +\Delta^H}{2}\right\|}_F^2+
\frac{1}{|\lambda|^2}{\left\|\frac{\Delta -\Delta^H}{2}\right\|}_F^2\,\right|~ \Delta \in \C^{n , n},\,
  \Delta x_2=r,~\Delta ^H x_1=s\right \}\\
&& =
\inf\left\{\left. {\|\Delta_R\|}_F^2+{\|\Delta_E\|}_F^2 ~\right|
  \Delta \in \C^{n , n},\, \Delta x_2=r,~\Delta ^H x_1=s,~\textstyle\Delta_R=-\frac{\Delta+\Delta^H}{2},~\Delta_E=\frac{\overline \lambda}{|\lambda|^2}\cdot\big(\frac{ \Delta-\Delta^H}{2}\big)\right \}
 \end{eqnarray*}
and hence
\begin{eqnarray}\label{equi_6}
  \inf\left\{ {\|\Delta_R\|}_F^2+{\|\Delta_E\|}_F^2 ~\big|~ \Delta_R,\,\Delta_E \in {\rm Herm}(n)
  ~{\rm satisfying}~\eqref{REequi_1}\,{\rm and}~\eqref{REequi_2}  \right \}\nonumber \\
\leq   \inf \left\{{\left\|\frac{\Delta +\Delta^H}{2}\right\|}_F^2+
\frac{1}{|\lambda|^2}{\left\|\frac{\Delta -\Delta^H}{2}\right\|}_F^2~\Big |~ \Delta \in \C^{n , n},\,
  \Delta u=r,~\Delta ^H w=s\right \}
 \end{eqnarray}
which finishes the proof of~\eqref{REequi_4}.

 \noindent
``$2)\Leftrightarrow 3)$'': By Theorem~\ref{thm:spcl_matrix_map}, there exist
 $\Delta \in \C^{n , n}$ satisfying $\Delta x_2=r$ and $\Delta^H x_1=s $
 if and only if $x_2^Hs=r^Hx_1$ which in turn is equivalent to $x_3^HB^Hx_1=0$.
\eproof

In contrast to Theorem~\ref{thm:strJE} and~\ref{thm:ustrRE}, we only obtain bounds
for the symmetry-structure-preserving eigenpair backward error $\eta^{{\mathcal S}}(R,E,\lambda,x)$.

\begin{theorem}\label{thm:strRE}
 Let $L(z)$ be a pencil as in~\eqref{eq:defL_cz}, and let
$\lambda \in i \R$ and $x \in \C^{2n+m}\setminus \{0\}$. Partition $x=[x_1^T~x_2^T~x_3^T]^T$
so that $x_1,\,x_2\in \C^n$ and $x_3 \in \C^{m}$, and set
$r=(J-R+\lambda E)x_2 +Bx_3$ and $s=-(J+R+\lambda E)x_1$.
Then $\eta^{{\mathcal S}}(R,E,\lambda,x)$ is finite if and only if $x_3=0$ and $B^Hx_1=0$.
If the latter conditions are satisfied then
\begin{eqnarray}\label{eq:RE_strberrbound1}
{\|\widehat{\Delta}\|}_F~\leq ~ \eta^{{\mathcal S}}(R,E,\lambda,x)~ \leq ~
\sqrt{{\left\|\frac{\widehat{\Delta} +\widehat{\Delta}^H}{2}\right\|}_F^2+
\frac{1}{|\lambda|^2}{\left\|\frac{\widehat{\Delta} -\widehat{\Delta}^H}{2}\right\|}_F^2}
\quad \quad {\rm if }~|\lambda|\leq 1
\end{eqnarray}
and
\begin{eqnarray}\label{eq:RE_strberrbound2}
\frac{{\|\widehat {\Delta}\|}_F}{|\lambda|}~\leq ~ \eta^{\mathcal S}(R,E,\lambda,x)~ \leq ~
\sqrt{{\left\|\frac{\widehat{\Delta} +\widehat{\Delta}^H}{2}\right\|}_F^2+
\frac{1}{|\lambda|^2}{\left\|\frac{\widehat{\Delta} -\widehat{\Delta}^H}{2}\right\|}_F^2}
\quad \quad {\rm if }~|\lambda|\geq 1,
\end{eqnarray}
where $\widehat{\Delta}$ is given by
\begin{equation*}
\widehat{\Delta}= \left\{\begin{array}{ll}
\frac{rx_2^H}{\|x_2\|^2} & \mbox{ if } x_1 = 0, \\[1ex]
\frac{x_1s^H}{\|x_1\|^2} & \mbox{ if } x_2 = 0,\\[1ex]
\frac{rx_2^H}{\|x_2\|^2} + \frac{x_1s^H}{\|x_1\|^2}\left(I_{n} - \frac{x_2x_2^H}{\|x_2\|^2}\right)
& \mbox{ otherwise. }  \end{array}\right.
\end{equation*}
\end{theorem}

\proof
Combining Remark~\ref{rem1:15.9.17} and Lemma~\ref{lem1:str_RE}, we obtain that
$\eta^{\mathcal S}(R,E,\lambda,x)$ is finite if and only if $x_3^HB^Hx_1=0$ and
$B^Hx_1+Sx_3=0$. The latter conditions hold if and only if $x_3=0$ and $B^Hx_1=0$, because
$S$ is definite. Thus let $x=[x_1^T~x_2^T~x_3^T]^T$ be such that $x_3=0$ and $B^Hx_1=0$.
Then we obtain from~\eqref{eq:deff_strRE_eperror} and by using Lemma~\ref{lem1:str_RE} that
\begin{eqnarray}\label{eq:15.9.17}
(\eta^{\mathcal S}(R,E,\lambda,x))^2\!\!\nonumber
&=&\inf\Big\{{\|[\Delta_R~\Delta_E]\|}_{F}^2 \,\Big |~
\Delta_R,\Delta_E \in {\rm Herm}(n),~{\rm satisfying}~\eqref{REequi_1}~{\rm and}~\eqref{REequi_2} \Big\}\\
&=& \inf \Bigg\{{\left\|\frac{\Delta +\Delta^H}{2}\right\|}_F^2\!\!\!+
\frac{1}{|\lambda|^2}{\left\|\frac{\Delta -\Delta^H}{2}\right\|}_F^2\Bigg | \Delta \in \C^{n , n},\,
  \Delta x_2=r,\,\Delta^Hx_1=s\Bigg \},\qquad\quad
\end{eqnarray}
where the last equality is due to Lemma~\ref{lem1:str_RE}.
Note that for any $\Delta \in \C^{n , n}$,  the Hermitian and skew-Hermitian parts of $\Delta$ satisfy
$\|\Delta\|_F^2={\big\|\frac{\Delta +\Delta^H}{2}\big\|}_F^2+
{\big\|\frac{\Delta -\Delta^H}{2}\big\|}_F^2$.
This implies
\begin{eqnarray}\label{eq:RE_auxi_strberrbound1}
{\|{\Delta}\|}_F^2~\leq ~ {\left\|\frac{\Delta +\Delta^H}{2}\right\|}_F^2+
\frac{1}{|\lambda|^2}{\left\|\frac{\Delta -\Delta^H}{2}\right\|}_F^2~ \quad \quad {\rm if }~|\lambda|\leq 1
\end{eqnarray}
and
\begin{eqnarray}\label{eq:RE_auxi_strberrbound2}
\frac{{\| {\Delta}\|}_F^2}{|\lambda|^2}~\leq ~ {\Big\|\frac{\Delta +\Delta^H}{2}\Big\|}_F^2+
\frac{1}{|\lambda|^2}{\Big\|\frac{\Delta -\Delta^H}{2}\Big\|}_F^2 ~
\quad \quad {\rm if }~|\lambda|\geq 1
\end{eqnarray}
for all $\Delta \in \C^{n , n}$. Then taking the infimum over all
$\Delta$ satisfying $\Delta x_2=r$ and $\Delta^Hx_1=s$ in~\eqref{eq:RE_auxi_strberrbound1}
and~\eqref{eq:RE_auxi_strberrbound2}, and by using the minimal Frobenius norm mapping from
Theorem~\ref{thm:spcl_matrix_map} we obtain~\eqref{eq:RE_strberrbound1} and~\eqref{eq:RE_strberrbound2}.
\eproof

\begin{example}\rm
The reason why we only obtain bounds in Theorem~\ref{thm:strRE} is the fact that the
infimum in~\eqref{eq:15.9.17}
need not be attained by the matrix $\widehat\Delta$ from Theorem~\ref{thm:strRE}. As an
example, consider the pencil $L(z)$ as in~\eqref{eq:defL_cz} with
\[
J=\mat{cc}0&-1\\ 1&0\rix,\;R=\mat{cc}0&0\\ 0&1\rix,\;E=B=\mat{cc}0&0\\ 0&0\rix,\;\mbox{and}\;S=I_2
\]
and let $\lambda=\frac{1}{4}i$ and $x=\mat{cccccc}0&0&1&1&0&0\rix^T$, i.e., $x_1=x_3=0\in\mathbb C^2$
and $x_2=\mat{cc}1&1\rix^T$. We then obtain $s=-(J+R+\lambda E)x_1=0$ as well as
\[
r=(J-R-\lambda E)x_2+Bx_3=\mat{c}-1\\ 0\rix\quad\mbox{and}\quad
\widehat\Delta=\frac{rx_2^H}{\|x_2\|^2}=\mat{cc}-\frac{1}{2}&-\frac{1}{2}\\ 0&0\rix
\]
which by~\eqref{eq:RE_strberrbound1} gives the bounds
\[
\frac{1}{2}=\|\widehat\Delta\|_F\leq\eta^{\mathcal S}(R,E,\lambda,x)\leq
\sqrt{{\left\|\mat{cc}-\frac{1}{2}&-\frac{1}{4}\\ -\frac{1}{4}&0\rix\right\|}_F^2+
\frac{1}{|\lambda|^2}{\left\|\mat{cc}0&-\frac{1}{4}\\ \frac{1}{4}&0\rix\right\|}_F^2}
=\sqrt{2.375}.
\]
On the other hand, for the Hermitian matrix
\[
\Delta:=\mat{cc}-1&0\\ 0&0\rix,
\]
we have $\Delta x_2=r$ and thus we obtain from~\eqref{eq:15.9.17} that
$\eta^{\mathcal S}(R,E,\lambda,x)\leq\|\Delta\|_F=1$.
\end{example}

It remains an open problem to determine the exact value for $\eta^{\mathcal S}(R,E,\lambda,x)$
and for the same reason, also the computation of the eigenvalue backward error
$\eta^{\mathcal S}(R,E,\lambda)$ is a challenging problem.

\subsection{Perturbations only in $J$ and $R$}

Next, we consider perturbations that only effect the blocks $J$ and $R$ in a pencil $L(z)$ as
in~\eqref{eq:defL_cz}. Let $\lambda\in \C$ and $x \in \C^{2n+m}\setminus \{0\}$.
Then by the terminology outlined in Section~\ref{sec:section_def_epbr}, the block- and
symmetry-structure-preserving eigenpair backward errors $\eta^{\mathcal B}(J,R,\lambda,x)$
and $\eta^{\mathcal S}(J,R,\lambda,x)$ are defined by
\begin{equation}\label{eq:deff_RE_eperror}
\eta^{\mathcal B}(J,R,\lambda,x)=\inf\Big\{\|[\Delta_J~\Delta_R]\|_{F} \, \Big |\,
\big((M-\Delta_M)+\lambda(N-\Delta_N)\big)x=0,\,\Delta_M + z \Delta_N \in \mathcal B\Big\},
\end{equation}
and
\begin{equation}\label{eq:def_strRE_eperr}
\eta^{\mathcal S}(J,R,\lambda,x)=\inf\Big\{\|[\Delta_J~\Delta_R]\|_{F} \, \Big |\,
\big((M-\Delta_M)+\lambda(N-\Delta_N)\big)x=0,\,\Delta_M + z \Delta_N \in \mathcal S\Big\},
\end{equation}
respectively, where $\mathcal B$ is the set of all pencils
of the form $\Delta L(z)=\Delta_M + z \Delta_N$ with the block-structure
\[
\Delta_M= \mat{ccc} 0 & \Delta_J-\Delta_R & 0\\ {(\Delta_J-\Delta_R )}^H &0&0 \\0&0&0\rix, \ \Delta_N=0,
\]
and $\Delta_J,\Delta_R\in\mathbb C^{n,n}$, while $\mathcal S$ is the corresponding set
of pencils $\Delta L(z)=\Delta_M + z \Delta_N$ as in~\eqref{14.9.17:pencil} with
$\Delta_J\in\text{SHerm}(n)$ and $\Delta_R\in\text{Herm}(n)$.
If the perturbations are restricted to be real, then the above backward errors are denoted by
$\eta^{\mathcal B_\R}(J,R,\lambda,x)$ and $\eta^{{\mathcal S_\R}}(J,R,\lambda,x)$, respectively.
As usual, we first investigate conditions for the finiteness of $\eta^{\mathcal B}(J,R,\lambda,x)$.

\begin{remark}\label{rem2:15.9.17}\rm
Let $L(z)$ be a pencil as in~\eqref{eq:defL_cz}, and let $\lambda \in\mathbb C$ and $x=[x_1^T~x_2^T~x_3^T]^T$
be such that $x_1,\,x_2 \in \C^{n}$ and $x_3\in \C^m$. Then for any $\Delta_J,\,
\Delta_R \in \C^{n , n}$ and corresponding $\Delta L(z)=\Delta_M +z \Delta_N \in \mathcal B$, we have
$(L-\Delta L)(\lambda)x=0$ if and only if
\begin{eqnarray}
(\Delta_J - \Delta_R)x_2&=&(J-R+\lambda E)x_2 +B x_3, \label{JRequi_1}\\
(\Delta_J- \Delta_R)^Hx_1 &=& (-J-R-\lambda E)x_1,\label{JRequi_2}\\
0&=& B^Hx_1+Sx_3. \label{JRequi_3}
\end{eqnarray}
Consequently, $\eta^{\mathcal B}(J,R,\lambda,x)$ is finite if and only if~\eqref{JRequi_1}--\eqref{JRequi_3}
are satisfied.
\end{remark}
\begin{lemma}\label{lem:JRstr_unstr}
Let $L(z)$ be a pencil as in~\eqref{eq:defL_cz}, and let
$\lambda \in i \R$ and $x \in \C^{2n+m}\setminus \{0\}$. Partition $x=[x_1^T~x_2^T~x_3^T]^T$
such that $x_1,\,x_2 \in \C^{n}$ and $x_3 \in \C^{m}$ and set $r=(J-R+\lambda E)x_2 +B x_3$
and $s=(-J-R-\lambda E)x_1$. Then the following statements are equivalent.
\begin{enumerate}
  \item There exist $\Delta_J,\, \Delta_R \in \C^{n , n}$ satisfying~\eqref{JRequi_1} and~\eqref{JRequi_2}.
  \item There exist $\Delta \in \C^{n , n}$ such that $\Delta x_2=r$ and $\Delta ^H x_1=s$.
  \item There exist $\Delta_J\in {\rm SHerm}(n),\, \Delta_R  \in {\rm Herm}(n)$
  satisfying~\eqref{JRequi_1} and~\eqref{JRequi_2}.
    \item The identity $x_3^HB^Hx_1=0$ is satisfied.
\end{enumerate}
Moreover,
\begin{eqnarray}\label{JRequi_4}
    \inf\left\{\left. {\|\Delta_J\|}_F^2+{\|\Delta_R\|}_F^2 ~\right|~ \Delta_J,\,\Delta_R \in \C^{n , n}~{\rm satisfy}~\eqref{JRequi_1}\,{\rm and}~\eqref{JRequi_2}  \right \}\nonumber\\
    =\inf \left\{\left.\frac{{\|\Delta\|}_F^2}{2}~\right|~ \Delta \in \C^{n , n},\,
  \Delta x_2=r,~\Delta ^H x_1=s\right \},
 \end{eqnarray}
and
\begin{eqnarray}\label{JRequi_5}
    \inf\Big\{ {\|\Delta_J\|}_F^2+{\|\Delta_R\|}_F^2 ~\Big|& \Delta_J\in {\rm SHerm}(n),\,
    \Delta_R \in {\rm Herm}(n)~{\rm satisfying}~\eqref{JRequi_1}\,{\rm and}~\eqref{JRequi_2}\Big \}\nonumber\\
&=\inf \left\{{\|\Delta\|}_F^2~\Big |~ \Delta \in \C^{n , n},\,
  \Delta x_2=r,~\Delta ^H x_1=s\right \}.
\end{eqnarray}
\end{lemma}

\proof
``$1)\Rightarrow 2)$'': Let $\Delta_J,\, \Delta_R \in \C^{n , n}$ be such that they satisfy~\eqref{JRequi_1}
and~\eqref{JRequi_2}. By setting $\Delta=\Delta_J-\Delta_R$ we get $\Delta x_2=r$ and ${\Delta}^H x_1=s$.
Furthermore, we have
\[
\|\Delta\|_F^2 \leq \big(\|\Delta_J \|_F+\|\Delta_R\|_F\big)^2
\leq 2\big({\|\Delta_J \|}_F^2+{\|\Delta_R\|}_F^2\big),
\]
where the last inequality is an elementary application of the Cauchy Schwartz inequality (in $\mathbb R^2)$.
But then the inequality ``$\geq$'' in~\eqref{JRequi_4} can be easily shown by following the arguments in the proof
of ``$1)\Rightarrow 2)$'' in Lemma~\ref{thm;equivalentcond}.

\noindent
``$2)\Rightarrow 1)$'': Suppose that $\Delta \in \C^{n, n}$ is such that
$\Delta x_2=r$ and ${\Delta}^Hx_1=s$ and define $\Delta_J=\frac{1}{2}\Delta$
and $\Delta_R=-\frac{1}{2}\Delta$. Then $\Delta_J$ and $\Delta_R$ satisfy~\eqref{JRequi_1} and~\eqref{JRequi_2}.
Also, we obtain $${\|\Delta_J\|}_F^2+{\|\Delta_R\|}_F^2=\frac{{\|\Delta\|}_F^2}{2}$$ and hence
``$\leq$'' in~\eqref{JRequi_4} can be easily shown by following the arguments of the proof of
``$2)\Rightarrow 1)$'' in Lemma~\ref{thm;equivalentcond}.

\noindent
``$2)\Rightarrow 3)$'': Let $\Delta \in \C^{n, n}$ be such that
$\Delta u=r$ and ${\Delta}^H w=s$. Then by setting $\Delta_J=\frac{\Delta-{\Delta}^H}{2}$
and $\Delta_R=-\frac{\Delta+{\Delta}^H}{2}$, we get $\Delta_J \in {\rm SHerm}(n)$,
$\Delta_R \in {\rm Herm}(n)$  such that~\eqref{JRequi_1} and~\eqref{JRequi_2} hold.
Furthermore, we have
\[
 {\|\Delta_J\|}_F^2+{\|\Delta_R\|}_F^2={\left\|\frac{\Delta-{\Delta}^H}{2}\right\|}_F^2+
 {\left\|\frac{\Delta+{\Delta}^H}{2}\right\|}_F^2= {\|\Delta\|}_F^2.
 \]
Thus, arguments similar to those in the proof of ``$2)\Rightarrow 1)$'' in Lemma~\ref{thm;equivalentcond}
give ``$\leq$'' in~\eqref{JRequi_5}.

\noindent
``$3)\Rightarrow 2)$'': Let $\Delta_J\in {\rm SHerm}(n)$ and $ \Delta_R \in {\rm Herm}(n)$ be such that they
satisfy~\eqref{JRequi_1} and~\eqref{JRequi_2}. Define $\Delta=\Delta_J-\Delta_R$
then $\Delta x_2=r$ and ${\Delta}^H x_1=s$. Note that $\Delta_J$ and $-\Delta_R$ are, respectively,
the unique skew-Hermitian and Hermitian parts of
$\Delta$, i.~e., $\Delta_J=\frac{\Delta-{\Delta}^H}{2}$ and $\Delta_R=-\frac{\Delta+{\Delta}^H}{2}$.
This implies
\[
 {\|\Delta_J\|}_F^2+{\|\Delta_R\|}_F^2={\left\|\frac{\Delta-{\Delta}^H}{2}\right\|}_F^2+
 {\left\|\frac{\Delta+{\Delta}^H}{2}\right\|}_F^2= {\|\Delta\|}_F^2.
 \]
Then again arguments similar to those in the proof of ``$1)\Rightarrow 2)$'' in Lemma~\ref{thm;equivalentcond}
give ``$\geq$'' in~\eqref{JRequi_5}.

\noindent
``$2)\Leftrightarrow 4)$'': This follows immediately from Theorem~\ref{thm:spcl_matrix_map}.
\eproof
The following theorem yields the values of $\eta^{\mathcal B}(J,R,\lambda,x)$, $\eta^{{\mathcal S}}(J,R,\lambda,x)$, and also
of their real counterparts if $L(z)$ is real. It also gives
the values of $\eta^{\mathcal B}(J,R,\lambda):=\inf_{x\in \C^{2n+m}\setminus\{0\}}\eta^{\mathcal B}(J,R,\lambda,x)$ and
$\eta^{\mathcal S}(J,R,\lambda):=\inf_{x\in \C^{2n+m}\setminus\{0\}}\eta^{\mathcal S}(J,R,\lambda,x)$.

\begin{theorem}\label{thm:ustrJR}
Let $L(z)$ be a pencil defined by~\eqref{eq:defL_cz},
$\lambda \in i \R$ and $x \in \C^{2n+m}\setminus \{0\}$. Partition $x=[x_1^T~x_2^T~x_3^T]^T$
such that $x_1,\,x_2 \in \C^n $ and $x_3 \in \C^{m}$, and set
$r=(J-R+\lambda E)x_2 +Bx_3$ and $s=-(J+R+\lambda E)x_1$. Then the following statements hold:
\begin{enumerate}
\item $\eta^{\mathcal B}(J,R,\lambda,x)$ and $\eta^{\mathcal S}(J,R,\lambda,x)$
are finite if and only if $x_3=0$ and $B^Hx_1=0$.
If the latter conditions are satisfied then
\begin{eqnarray*}\label{eq:JR_strberr}
\eta^{\mathcal B}(J,R,\lambda,x)~ = ~
\frac{{\|\widehat {\Delta}\|}_F}{\sqrt{2}} \quad \text{and}\quad
\eta^{\mathcal S}(J,R,\lambda,x)~ = ~{\|\widehat {\Delta}\|_F},
\end{eqnarray*}
as well as
\begin{equation*}
\eta^{\mathcal B}(J,R,\lambda)=\frac{\sigma_{\min}(J-R+\lambda E)}{\sqrt{2}} \quad \text{and}\quad
\eta^{\mathcal S}(J,R,\lambda)=\sigma_{\min}(J-R+\lambda E),
\end{equation*}
where $\widehat{\Delta}$ is given by
\begin{equation*}
\widehat{\Delta}= \left\{\begin{array}{ll}
\frac{rx_2^H}{\|x_2\|^2} & \mbox{ if } x_1 = 0, \\[1ex]
\frac{x_1s^H}{\|x_1\|^2} & \mbox{ if } x_2 = 0,\\[1ex]
\frac{rx_2^H}{\|x_2\|^2} + \frac{x_1s^H}{\|x_1\|^2}\left(I_{n} - \frac{x_2x_2^H}{\|x_2\|^2}\right)
& \mbox{ otherwise. }  \end{array}\right.
\end{equation*}
\item  Suppose that $L(z)$ is real. If ${\rm rank}\left([x_1~\overline x_1]\right)=
{\rm rank}\left([x_2~\overline x_2]\right)=2$ then
$\eta^{\mathcal B_\R}(J,R,\lambda,x)$ and $\eta^{{\mathcal S}_\R}(J,R,\lambda,x)$
are finite if and only if $x_3=0$, $B^Tx_1=0$ and $\lambda x_2^TEx_1=0$. If the latter conditions
are satisfied then
\begin{eqnarray}\label{eq:JRreal_strberr}
 \eta^{{\mathcal B}_\R}(J,R,\lambda,x) =\frac{{\|\widetilde {\Delta} \|}_F}{\sqrt{2}}
\quad {\rm and}\quad
\eta^{{{\mathcal S}_\R}}(J,R,\lambda,x) ={\|\widetilde {\Delta} \|}_F,
\end{eqnarray}
where $\widetilde{\Delta} \in \R^{n, n}$ is given by
\begin{equation*}\label{eq:JRreal_sepbrrobspecial}
\widetilde{\Delta}= [r~ \overline r][x_2~ \overline x_2]^{\dag}+
\big([s~ \overline s][x_1~ \overline x_1]^{\dag}\big)^H
-\big([s~ \overline s][x_1~ \overline x_1]^{\dag}\big)^H\big([x_2~ \overline x_2][x_2~ \overline x_2]^{\dag}\big).
\end{equation*}
\end{enumerate}
\end{theorem}

\proof The proof of $1)$ follows the same lines as that of Theorem~\ref{thm:ustrJE} by using
Lemma~\ref{lem:JRstr_unstr} and Theorem~\ref{thm:spcl_matrix_map}.

\noindent
Concerning the proof of 2), recall that when $L(z)$ is real, then $\eta^{\mathcal B_{\R}}(J,R,\lambda,x)$ is the
eigenpair backward error obtained by allowing only real perturbations to the blocks $J$ and $R$ of $L(z)$.
Now for any $\Delta_J,\, \Delta_R \in \R^{n , n}$ and corresponding real
$\Delta L(z)=\Delta_M +z \Delta_N \in \mathcal B$, we have $(L-\Delta L)(\lambda)x=0$ if and only if
\begin{eqnarray}
(\Delta_J - \Delta_R){x_2} &=&{(J-R+\lambda E)x_2 +B x_3}\label{realJRequi_1}\\
(\Delta_J- \Delta_R)^T {x_1} &=& {(-J-R-\lambda E)x_1}\label{realJRequi_2}\\
0&=& B^Tx_1+Sx_3. \label{realJRequi_3}
\end{eqnarray}
Since $\Delta_J$ and $\Delta_R$ are real,~\eqref{realJRequi_1} and~\eqref{realJRequi_2} can be
equivalently written as
\begin{eqnarray}\label{realJRequi_4}
(\Delta_J - \Delta_R)[x_2~\bar x_2] =[r~\bar r] \quad {\rm and} \quad (\Delta_J- \Delta_R)^T [x_1~\bar x_1]= [s~\bar s].
 \end{eqnarray}
Following the lines of the proof of Lemma~\ref{lem:JRstr_unstr}, there exist real matrices $\Delta_J$ and
$\Delta_R$ satisfying~\eqref{realJRequi_4} if and only if there exist $\Delta \in \R^{n , n}$ such that
$\Delta[x_2~\bar x_2] =[r~\bar r]$ and ${\Delta}^T[x_1~\bar x_1] =[s~\bar s]$. Applying
Theorem~\ref{cor:real_spcl_cor}, we find that this is the case if and only if
\[
x_2^Hs=r^Hx_1\quad\mbox{and}\quad x_2^Ts=r^Tx_1
\]
which, using the definition of $r$ and $s$, is in turn equivalent to the conditions
$$x_3^HB^Hx_1=0\quad\mbox{and}\quad2\lambda x_2^TEx_1=x_3^TB^Tx_1.$$ The latter conditions
together with $B^Tx_1+Sx_3=0$ give $x_3=0$, $B^Tx_1=0$ and $\lambda x_2^TEx_1=0$, because $S$ is assumed to be
positive definite. Therefore from~\eqref{eq:deff_RE_eperror}, $\eta^{\mathcal B_{\R}}(J,R,\lambda,x)$
is finite if and only if $x$ satisfies $x_3=0$, $B^Tx_1=0$ and $\lambda x_2^TEx_1=0$. If this is
the case, then we find that
\begin{eqnarray*}\label{realJRequi_5}
\eta^{\mathcal B_{\R}}(J,R,\lambda,x)&=&
\inf \Big\{ {\|[\Delta_J~\Delta_R]\|}_{F}\, \Big|\,\Delta_J,\,\Delta_R \in \R^{n , n} \,
{\rm satisfy}\,\eqref{realJRequi_4}\Big\}\\
&=&\inf \left\{ \left.\frac{{\|\Delta\|}_F}{\sqrt{2}} \right|\,\Delta \in \R^{n , n},~ \Delta [x_2 ~\bar x_2]=
[r ~\bar r]~{\rm and}~{\Delta}^T[x_1 ~\bar x_1]=[s ~\bar s]\right\}.
\end{eqnarray*}
Thus~\eqref{eq:JRreal_strberr} follows for $\eta^{\mathcal B_{\R}}(J,R,\lambda,x)$ by using Theorem~\ref{cor:real_spcl_cor}.
Similarly we can also establish~\eqref{eq:JRreal_strberr} for $\eta^{\mathcal S_{\R}}(J,R,\lambda,x)$.
\eproof

\subsection{Perturbation only to J and B, or  R and B, or E and B}
In this section, we obtain block-structure-preserving eigenpair or eigenvalue backward errors
when only the blocks $J$ and $B$ in a pencil $L(z)$ as in~\eqref{eq:defL_cz} are perturbed.
Unfortunately, is seems that this approach cannot be generalized to obtain the correpsonding
symmetry-structure-preserving backward errors.

Let $\lambda \in \C$ and $x \in \C^{2n+m}\setminus \{0\}$, then by the terminology outlined in
Section~\ref{sec:section_def_epbr}, the block-structure-preserving eigenpair backward error
$\eta^{\mathcal B}(J,B,\lambda,x)$ is defined by
\begin{eqnarray}\label{eq:def_JB_eperr}
\eta^{\mathcal B}(J,B,\lambda,x)=\inf\Big\{{\|[\Delta_J~\Delta B]\|}_{F}\, \Big | &
\Delta_J\in \C^{n , n},\,\Delta B \in \C^{n , m},\,\Delta_M + z \Delta_N \in \mathcal B, \nonumber\\
&\big((M-\Delta_M)+\lambda(N-\Delta_N)\big)x=0\Big\},
\end{eqnarray}
where $\mathcal B$ is the set of all pencils of the form $\Delta L(z)=
\Delta_M + z \Delta_N$ with $$\Delta_M= \mat{ccc} 0 & \Delta_J& \Delta_B\\ \Delta_J^H &0&0 \\\Delta_B^H&0&0\rix
\quad\mbox{and}\quad\Delta_N=0.$$
If the perturbations are restricted to be real then the above error is denoted by
$\eta^{\mathcal B_\R}(J,B,\lambda,x)$.

\begin{remark}\rm
Let $L(z)$ be a pencil as in~\eqref{eq:defL_cz}, and let $\lambda \in\mathbb C$ and $x=[x_1^T~x_2^T~x_3^T]^T$
be such that $x_1,\,x_2 \in \C^{n}$ and $x_3\in \C^m$.
Then for any $\Delta_J \in \C^{n , n}$, $\Delta_B \in \C^{n , m}$, and corresponding
$\Delta L(z)=\Delta_M +z \Delta_N \in \mathcal B$, we have $(L-\Delta L)(\lambda)x=0$ if and only if
\begin{eqnarray*}
\Delta_Jx_2 + \Delta_B x_3 &=&(J-R+\lambda E)x_2 +B x_3, \\
\Delta_J^H x_1 &=&{(-J-R-\lambda E)x_1},\\
\Delta_B^H x_1&=& B^Hx_1+Sx_3,
\end{eqnarray*}
which in turn is equivalent to
\begin{eqnarray}
 \mat{cc} \Delta_J &\Delta_B \rix \underbrace{\mat{c} x_2\\x_3 \rix}_{=u} &=&
 \underbrace{(J-R+\lambda E)x_2 +B x_3}_{=r}, \label{JBequi_1}\\
 \mat{cc} \Delta_J &\Delta_B \rix^H \underbrace{x_1}_{=w} &=&\underbrace{
 \mat{c} -(J+R+\lambda E)x_1\\ B^Hx_1 + Sx_3 \rix}_{=s}.\label{JBequi_2}
\end{eqnarray}
In particular, $\eta^{\mathcal B}(J,B,\lambda,x)$ is finite if and only if~\eqref{JBequi_1}--\eqref{JBequi_2}
are satisfied.
\end{remark}

\begin{lemma}\label{lem1:unstr_JB}
Let $L(z)$ be a pencil as in~\eqref{eq:defL_cz}, and let $\lambda \in i \R$ and $x \in \C^{2n+m}\setminus \{0\}$.
Partition $x=[x_1^T~x_2^T~x_3^T]^T$ such that $x_1,\,x_2 \in \C^{n}$ and $x_3 \in \C^{m}$, and let
$u\,,w\,,r$ and $s$ be defined as in~\eqref{JBequi_1} and~\eqref{JBequi_2}. Then the following
statements are equivalent.
\begin{enumerate}
  \item There exist $\Delta_J \in \C^{n,n}$ and $\Delta_B \in \C^{n , m}$ satisfying~\eqref{JBequi_1} and~\eqref{JBequi_2}.
  \item  There exist $\Delta \in \C^{n , n+m}$ such that $\Delta u=r$ and $\Delta^H w=s$.
  \item $x$ satisfies $x_3=0$.
\end{enumerate}
 Moreover, we have
\begin{eqnarray*}\label{JBequi_3}
\inf\Big\{ {\|\Delta_J\|}_F^2+{\|\Delta_B\|}_F^2 \,\Big| ~\Delta_J \in \C^{n,n},\,\Delta_E \in\C^{n , m}~
  {\rm satisfy}~\eqref{JBequi_1}~{\rm and}~\eqref{JBequi_2}  \Big \}\\
  =\inf \left\{\left.{\|\Delta\|}_F^2~\right |~ \Delta \in \C^{n ,n+m},\,\Delta u=r,~\Delta ^H w=s\right \}.
 \end{eqnarray*}
\end{lemma}
\proof ``$1)\Rightarrow 2)$'' is obvious while ``$2)\Rightarrow 3)$'' is implied by
Theorem~\ref{thm:spcl_matrix_map} using the fact that $S$ is definite. The last part then
follows from the observation that any $\Delta \in \C^{n , n+m}$ can be written as
$\Delta=[{\Delta}_1~~{\Delta}_2]$,
where ${\Delta}_1 \in \C^{n,n}$ and ${\Delta}_2\in \C^{n , m}$ such that
\[
{\|\Delta\|}_F={\|[{\Delta}_1\,~{\Delta}_2]\|}_F=\sqrt{{{\|{\Delta}_1\|}_F^2+{\|{\Delta}_2\|}}_F^2}.
\quad\mbox{\eproof}
\]
\begin{theorem}\label{thm:unstrerrorJB}
Let $L(z)$ be a pencil as in~\eqref{eq:defL_cz}, let
$\lambda \in i \R$ and $x \in \C^{2n+m}\setminus \{0\}$. Partition $x=[x_1^T~x_2^T~x_3^T]^T$
such that $x_1,\,x_2\in \C^n$ and $x_3 \in \C^{m}$ and set $ u=[x_2^T~x_3^T]^T$, $w=x_1$,
\begin{eqnarray*}
r=(J-R+\lambda E)x_2 +Bx_3\quad \text{and} \quad s=[-((J+R+\lambda E)x_1)^T~(B^Hx_1+Sx_3)^T]^T.
\end{eqnarray*}
Then the following statements hold.
\begin{enumerate}
\item $\eta^{\mathcal B}(J,B,\lambda,x)$ is finite if and only if $x_3=0$.
In that case, we have
\begin{eqnarray*}\label{eq:JB_unstrberr}
 \eta^{\mathcal B}(J,B,\lambda,x)~ = ~\sqrt{{ \|\widehat{ \Delta}_1\|}_F^2+{\|\widehat{\Delta}_2\|}_F^2},
\end{eqnarray*}
and
\begin{equation*}
 \eta^{\mathcal B}(J,B,\lambda)=\min\left\{\sigma_{\min}\left(\mat{cc}J-R+\lambda E&  B\rix^H\right),
 \sigma_{\min}(J-R+\lambda E)\right\},
\end{equation*}
where $\widehat{\Delta}_1$ and $\widehat {\Delta}_2$ are given by
\begin{equation*}
[\widehat{\Delta}_1~\widehat{\Delta}_2]= \left\{\begin{array}{ll}
\frac{ru^H}{\|u\|^2} & \mbox{ if } x_1 = 0, \\[1ex]
\frac{ws^H}{\|w\|^2} & \mbox{ if } x_2 = 0,\\[1ex]
\frac{ru^H}{\|u\|^2} + \frac{ws^H}{\|w\|^2}\left(I_{n+m} - \frac{uu^H}{\|u\|^2}\right)
& \mbox{ otherwise. }  \end{array}\right.
\end{equation*}
\item Suppose that $L(z)$ is real. If $\operatorname{rank}\left([x_1~\overline x_1]\right)
=\operatorname{rank}\left([x_2~\overline x_2]\right)=2$ then
$\eta^{{\mathcal B}_\R}(J,B,\lambda,x)$ is finite if and only if $x_3=0$ and $\lambda x_2^TEx_1=0$.
If the latter conditions are satisfied, then
\begin{eqnarray}\label{eq:JBreal_strberr}
 \eta^{{\mathcal B}_\R}(J,B,\lambda,x) =\sqrt{{\|\wt {\Delta}_1 \|}_F^2+{\|\wt {\Delta}_2 \|}_F^2},
\end{eqnarray}
where $\widetilde {\Delta}_1 \in \R^{n,n}$ and $\widetilde {\Delta}_2 \in \R^{n, m}$ are given by
\begin{equation*}
[\widetilde {\Delta}_1~\widetilde {\Delta}_2]= [r~ \bar r][u~ \bar u]^{\dag}+
\big([s~ \bar s][w~ \bar w]^{\dag}\big)^H
-\big([s~ \bar s][w~ \bar w]^{\dag}\big)^H\big([u~ \bar u][u~ \bar u]^{\dag}\big).
\end{equation*}
\end{enumerate}
\end{theorem}

\proof
The proof is analogous to the one of Theorem~\ref{thm:ustrJE} by using Lemma~\ref{lem1:unstr_JB} as
well as Theorem~\ref{thm:spcl_matrix_map} in the complex case and Theorem~\ref{cor:real_spcl_cor}
in the real case.
\eproof

\begin{remark}\label{thm:ustrRB}
{\rm
A result similar to Theorem~\ref{thm:unstrerrorJB}  can be obtained
for the complex and real block-structure-preserving eigenpair backward errors $\eta^{\mathcal B}(R,B,\lambda,x)$
and $\eta^{\mathcal B_\R}(R,B,\lambda,x)$ of a pair $(\lambda,x)\in (i\R)\times( \C^{2n+m}\setminus \{0\})$
when only the blocks $R$ and $B$ in a pencil $L(z)$ as in~\eqref{eq:defL_cz} are subject to perturbation. In fact, one easily obtains
 \[
  \eta^{\mathcal B}(R,B,\lambda,x)=\eta^{\mathcal B}(J,B,\lambda,x)\quad {\rm and} \quad
  \eta^{{\mathcal B}_\R}(R,B,\lambda,x)=\eta^{{\mathcal B}_\R}(J,B,\lambda,x).
 \]
 As a consequence we also have
 \[
 \eta^{\mathcal B}(R,B,\lambda)=\eta^{\mathcal B}(J,B,\lambda).
 \]
 }
\end{remark}

Finally, also the backward errors $\eta^{\mathcal B}(E,B,\lambda,x)$ and $\eta^{\mathcal B}(E,B,\lambda)$
with respect to perturbations only in the blocks $E$ and $B$ of $L(z)$ as in~\eqref{eq:defL_cz}
can be obtained in a similar manner. Since the actual result differs slightly from the previous
formulas, we present it as a theorem, but we omit the proof, since it is similar to the one of
Theorem~\ref{thm:ustrJE}.

\begin{theorem}\label{thm:ustrEB}
Let $L(z)$ be a pencil as in~\eqref{eq:defL_cz}, and let
$\lambda \in i \R$ and $x \in \C^{2n+m}\setminus \{0\}$. Partition $x=[x_1^T~x_2^T~x_3^T]^T$
such that $x_1,\,x_2 \in \C^n $ and $x_3 \in \C^{m}$, and set $u=[\lambda x_2^T~x_3^T]^T$, $w=x_1$,
\begin{eqnarray*}
  r=(J-R+\lambda E)x_2 +Bx_3\quad \text{and} \quad s=\Big[\frac{1}{\lambda}((J+R+\lambda E)x_1)^T~(B^Hx_1+Sx_3)^T\Big]^T.
\end{eqnarray*}
Then the following statements hold.
\begin{enumerate}
 \item $\eta^{\mathcal B}(E,B,\lambda,x)$ is finite if and only if $x_3=0$.
In that case, we have
\begin{eqnarray*}\label{eq:EB_strberr}
 \eta^{\mathcal B}(E,B,\lambda,x)~ = ~\sqrt{{ \|\widehat{ \Delta}_1\|}_F^2+{\|\widehat{\Delta}_2\|}_F^2},
\end{eqnarray*}
and
\begin{equation*}
 \eta^{\mathcal B}(E,B,\lambda)=\min\left\{\sigma_{\min}\left(\mat{cc}-\frac{(J-R+\lambda E)}{\lambda}&  B\rix^H\right),
 \frac{\sigma_{\min}(J-R+\lambda E)}{|\lambda|}\right\},
\end{equation*}
where $\widehat{\Delta}_1\in\mathbb C^{n,n}$ and $\widehat {\Delta}_2\in\mathbb C^{n,m}$ are given by
\begin{equation*}
[\widehat{\Delta}_1~\widehat{\Delta}_2]= \left\{\begin{array}{ll}
\frac{ru^H}{\|u\|^2} & \mbox{ if } x_1 = 0, \\[1ex]
\frac{ws^H}{\|w\|^2} & \mbox{ if } x_2 = 0,\\[1ex]
\frac{ru^H}{\|u\|^2} + \frac{ws^H}{\|w\|^2}\left(I_{n+m} - \frac{uu^H}{\|u\|^2}\right)
& \mbox{ otherwise. }  \end{array}\right.
\end{equation*}
\item Suppose that $L(z)$ is real. If $\operatorname{rank}\left([x_1~\overline x_1]\right)=
\operatorname{rank}\left([x_2~\overline x_2]\right)=2$ then
$\eta^{{\mathcal B}_\R}(E,B,\lambda,x)$ is finite if and only if $x_3=0$ and $\lambda x_2^TEx_1=0$.
If the latter conditions are satisfied then
\begin{eqnarray*}\label{eq:EBreal_strberr}
 \eta^{{\mathcal B}_\R}(E,B,\lambda,x) =\sqrt{{\|\widetilde {\Delta}_1 \|}_F^2+{\|\widetilde {\Delta}_2 \|}_F^2},
\end{eqnarray*}
where $\widetilde {\Delta}_1\in \R^{n,n}$ and $\widetilde {\Delta}_2 \in \R^{n, m}$ are given by
\begin{equation*}
[\widetilde {\Delta}_1~\widetilde {\Delta}_2]= [r~ \overline r][u~ \overline u]^{\dag}+
\big([s~ \overline s][w~ \overline w]^{\dag}\big)^H
-\big([s~ \overline s][w~ \overline w]^{\dag}\big)^H\big([u~ \overline u][u~ \overline u]^{\dag}\big).
\end{equation*}
\end{enumerate}
\end{theorem}

\section{Perturbation in any three of the matrices $J$, $R$, $E$ and $B$}
In this section, we define and compute block- and
symmetry-structure-preserving eigenpair or eigenvalue backward errors for pencils $L(z)$ as in~\eqref{eq:defL_cz},
while we consider perturbations in any three of the blocks $J,R,E,B$ of $L(z)$.

\subsection{Perturbations in the blocks J, R, and B}

We first concentrate on the case that perturbations are allowed to affect only the blocks
$J$, $R$, and $B$ of  a pencil $L(z)$ as in~\eqref{eq:defL_cz}.
If $\lambda \in \C$ and $x \in \C^{2n+m}\setminus\{0\}$, then following the terminology of
Section~\ref{sec:section_def_epbr}, the block- and symmetry-structure-preserving eigenpair backward errors
$\eta^{\mathcal B}(J,R,B,\lambda,x)$ and $\eta^{\mathcal S}(J,R,B,\lambda,x)$, respectively, are defined by
\begin{eqnarray*}\label{eq:def_JRB_eperr}
\eta^{\mathcal B}(J,R,B,\lambda,x)&=&\inf\Big\{{\big\|[\Delta_J~\Delta_R~\Delta_B]\big\|}_{F}\,\Big |\,
\big((M-\Delta_M)+\lambda(N-\Delta_N)\big)x=0,\,\Delta_M + z \Delta_N \in \mathcal B\Big\},\\
\eta^{\mathcal S}(J,R,B,\lambda,x)&=&\inf\Big\{{\big\|[\Delta_J~\Delta_R~\Delta_B]\big\|}_{F}\,\Big |\,
\big((M-\Delta_M)+\lambda(N-\Delta_N)\big)x=0,\,\Delta_M + z \Delta_N \in \mathcal S\Big\},
\end{eqnarray*}
where $\mathcal B$ denotes the set of all pencils of the form $\Delta L(z)=
\Delta_M + z \Delta_N$ with
\[
\Delta_M= \mat{ccc} 0 & \Delta_J-\Delta_R & \Delta_B\\ {(\Delta_J-\Delta_R )}^H &0&0 \\(\Delta_B)^H&0&0\rix,\ \Delta_N=0,
\]
and $\Delta_J,\Delta_R\in\mathbb C^{n,n}$, $\Delta_B\in\mathbb C^{n,m}$, while $\mathcal S$ denotes
the corresponding set of pencils that satisfy in addition $\Delta_J\in\text{SHerm}(n)$
and $\Delta_R\in\text{Herm}(n)$.
If the perturbations are restricted to be real then the above backward errors are denoted by
$\eta^{\mathcal B_\R}(J,R,B,\lambda,x)$ and $\eta^{{\mathcal S_\R}}(J,R,B,\lambda,x)$, respectively.

\begin{remark}
Let $L(z)$ be a pencil as in~\eqref{eq:defL_cz}, and let $\lambda\in\mathbb C$ and
$x=[x_1^T~x_2^T~x_3^T]^T$ such that $x_1,x_2 \in \C^n$ and $x_3 \in \C^m$.
Then for any $\Delta_J,\, \Delta_R \in \C^{n,n}$ and $\Delta_B \in \C^{n , m}$, and
corresponding $\Delta L(z)=\Delta_M +z \Delta_N \in \mathcal B$, we have
$(L-\Delta L)(\lambda)x=0$ if and only if
\begin{eqnarray*}
(\Delta_J-\Delta_R) x_2 + \Delta_B x_3 &=&(J-R+\lambda E)x_2 +B x_3, \\
(\Delta_J-\Delta_R)^Hx_1 &=&{(-J-R-\lambda E)x_1},\\
(\Delta_B)^Hx_1&=& B^Hx_1+Sx_3,
 \end{eqnarray*}
if and only if
\begin{eqnarray}
 \mat{cc} \Delta_J-\Delta_R &\Delta_B \rix \underbrace{\mat{c} x_2\\x_3 \rix}_{=u} &=&
 \underbrace{(J-R+\lambda E)x_2 +B x_3}_{=r} \label{JRBequi_1}, \\
 \mat{cc} \Delta_J-\Delta_R &\Delta_B \rix^H \underbrace{x_1}_{=w} &=&
 \underbrace{\mat{c} -(J+R+\lambda E)x_1\\ B^Hx_1 + Sx_3 \rix}_{=s}. \label{JRBequi_2}
\end{eqnarray}
\end{remark}

\begin{lemma}\label{lem:JRBbstr_unstr}
   Let $L(z)$ be a pencil defined by~\eqref{eq:defL_cz},
   $\lambda \in i \R$ and $x \in \C^{2n+m}\setminus \{0\}$. Partition $x=[x_1^T~x_2^T~x_3^T]^T$
 such that $x_1,\,x_2\in \C^n$ and $x_3 \in \C^{m}$ and let $u\,,w\,,r$ and $s$ be as defined in~\eqref{JRBequi_1}
 and~\eqref{JRBequi_2}. Then the following statements are equivalent.
\begin{enumerate}
  \item There exist $\Delta_J,\, \Delta_R\in \C^{n,n}$ and $\Delta_B \in \C^{n , m}$
  satisfying~\eqref{JRBequi_1} and~\eqref{JRBequi_2}.
  \item There exists ${\Delta} \in \C^{n , n+m}$ such that ${\Delta} u=r$ and $\Delta^H w=s$.
  \item There exist $\Delta_J\in {\rm SHerm}(n),\, \Delta_R  \in {\rm Herm}(n)$ and $\Delta_B \in \C^{n , m}$
  satisfying~\eqref{JRBequi_1} and~\eqref{JRBequi_2}.
  \item $x$ satisfies $x_3=0$.
\end{enumerate}
Moreover, we have
\begin{eqnarray*}\label{JRBequi_3}
    \inf\left\{ {\|[\Delta_J~\Delta_R~\Delta_B]\|}_F^2 \,\Big|\,
    \Delta_J,\,\Delta_R \in \C^{n,n},\, \Delta_B \in \C^{n , m}~{\rm satisfy}~\eqref{JRBequi_1}~
    {\rm and}~\eqref{JRBequi_2}  \right \}\qquad\nonumber\\
    =\inf \left\{\left.{{\frac{{\| {\Delta}_1\|}_F^2}{2}+{\|{\Delta}_2\|}_F^2}}~\right|~
    {\Delta}_1\in \C^{n,n},\,{\Delta}_2 \in \C^{n , m},\,  [{\Delta}_1~~{\Delta}_2] u=r,
    ~{[{\Delta}_1~{\Delta}_2]}^H w=s\right \},
\end{eqnarray*}
and
\begin{eqnarray*}\label{JRBequi_4}
&\inf\left\{ {\|[\Delta_J~\Delta_R~\Delta_B]\|}_F^2 \, \Big|\, \Delta_J\in {\rm SHerm}(n),
\,\Delta_R \in {\rm Herm}(n),\,\Delta_B \in \C^{n , m}\,{\rm satisfy}\,\eqref{JRBequi_1}\,{\rm and}
\,\eqref{JRBequi_2}  \right \}\nonumber\\
&=\inf \left\{{{{{\| {\Delta}_1\|}_F^2}+{\|{\Delta}_2\|}_F^2}}\,\Big |~
{\Delta}_1\in \C^{n,n},\,{\Delta}_2 \in \C^{n , m},\,  [{\Delta}_1~~{\Delta}_2] u=r,
    ~{[{\Delta}_1~{\Delta}_2]}^H w=s\right \}.
\end{eqnarray*}
\end{lemma}

\proof As seen in the proof of Lemma~\ref{lem1:unstr_JB}
any $\Delta \in \C^{n , n+m}$ can be written as
$\Delta=[{\Delta}_1~\,{\Delta}_2]$ where ${\Delta}_1 \in \C^{n,n}$ and ${\Delta}_2\in \C^{n , m}$
such that ${\|\Delta\|}_F={\|[{\Delta}_1\,~{\Delta}_2]\|}_F
=\sqrt{{{\|{\Delta}_1\|}_F^2+{\|{\Delta}_2\|}}_F^2}$.
With this key observation the proof is obtained by following exactly the same arguments
as in the proof of Lemma~\ref{lem:JRstr_unstr}.
\eproof

\begin{theorem}\label{thm:ustrJRB}
Let $L(z)$ be a pencil as in~\eqref{eq:defL_cz}, and let
$\lambda \in i \R$ and $x \in \C^{2n+m}\setminus \{0\}$. Partition $x=[x_1^T~x_2^T~x_3^T]^T$
such that $x_1,\,x_2 \in \C^n$ and $x_3 \in \C^{m}$ and define $ \hat u=[\sqrt{2} x_2^T~\,x_3^T]^T,~w=x_1$,
\begin{eqnarray*}
r=(J-R+\lambda E)x_2 +Bx_3\quad \text{and}\quad
\hat s=\big[-\frac{1}{\sqrt{2}}((J+R+\lambda E)x_1)^T~~(B^Hx_1+Sx_3)^T\big]^T.
  \end{eqnarray*}
Then the following statements hold:
\begin{enumerate}
\item $\eta^{\mathcal B}(J,R,B,\lambda,x)$ is finite if and only if $x_3=0$. In that case, we have
\begin{eqnarray*}\label{eq:JRB_unstrberr}
 \eta^{\mathcal B}(J,R,B,\lambda,x) =
{\sqrt{{\|\widehat {\Delta}_1\|}_F^2+{\|\widehat{\Delta}_2\|}_F^2}},
\end{eqnarray*}
and
\begin{equation*}
 \eta^{\mathcal B}(J,R,B,\lambda)=\min\left\{\sigma_{\min}\left(\mat{cc}\frac{(J-R+\lambda E)}{\sqrt{2}}&  B\rix^H\right),
 \frac{\sigma_{\min}(J-R+\lambda E)}{\sqrt{2}}\right\},
\end{equation*}
where $\widehat{\Delta}_1$ and $\widehat{\Delta}_2$ are given by
\begin{equation*}\label{eq:unstr_JRBsepbrrobspecial}
[\widehat{\Delta}_1~\widehat{\Delta}_2]= \left\{\begin{array}{ll}
\frac{r\hat u^H}{\|\hat u\|^2} & \mbox{ if } x_1 = 0, \\[1ex]
\frac{w\hat s^H}{\|w\|^2} & \mbox{ if } x_2 = 0,\\[1ex]
\frac{r\hat u^H}{\|\hat u\|^2} + \frac{w\hat s^H}{\|w\|^2}\left(I_{n+m} - \frac{\hat u \hat u^H}{\|\hat u\|^2}\right)
& \mbox{ otherwise. }  \end{array}\right.
\end{equation*}
\item If  $L(z)$ is real, and $\operatorname{ rank}\left([x_1~\overline x_1]\right)=
\operatorname{rank}\left([x_2~\overline x_2]\right)=2$, then $\eta^{{\mathcal B}_\R}(J,R,B,\lambda,x)$
is finite if and only if $x_3=0$ and $\lambda x_2^TEx_1=0$. If the latter conditions are satisfied then
\begin{eqnarray*}\label{eq:JRBreal_unstrberr}
 \eta^{{\mathcal B}_\R}(J,R,B,\lambda,x) =
{\sqrt{{\|\widetilde {\Delta}_1\|}_F^2+{\|\widetilde{\Delta}_2\|}_F^2}},
\end{eqnarray*}
where $\widetilde{\Delta}_1 \in \R^{n,n}$ and $\widetilde{\Delta}_2 \in \R^{n , m}$ are given by
\begin{equation*}\label{eq:unstr_JRBreal_sepbrrobspecial}
[\widetilde{\Delta}_1~\widetilde{\Delta}_2]= [r~ \overline r][\hat u~ \overline {\hat u}]^{\dag}+
\big([\hat s~ \overline {\hat s}][w~ \overline w]^{\dag}\big)^H
-\big([\hat s~ \overline {\hat s}][w~ \overline w]^{\dag}\big)^H\big([\hat u~ \overline {\hat u}][\hat u~ \overline {\hat u}]^{\dag}\big).
\end{equation*}
\end{enumerate}
\end{theorem}

\proof
Observe that if $u=[x_2^T~~x_3^T]^T$ and $s=\big[-((J+R+\lambda E)x_1)^T\,~(B^Hx_1+Sx_3)^T\big]^T$, then
\begin{eqnarray*}
\inf \left\{\left.{{\frac{{\| {\Delta}_1\|}_F^2}{2}+{\|{\Delta}_2\|}_F^2}}~\right|~
    {\Delta}_1,\,{\Delta}_2 \in \C^{n , n},\,  [{\Delta}_1~~{\Delta}_2] u=r,
    ~{[{\Delta}_1~{\Delta}_2]}^H w=s\right \}\\
   =\inf \left\{{{{\| { {\hat \Delta}}_1\|}_F^2+{\|{ {\hat \Delta}}_2\|}_F^2}}~\Big |~
    { {\hat \Delta}}_1,\,{ {\hat \Delta}}_2 \in \C^{n , n},\,  [{ {\hat \Delta}}_1~~{ {\hat \Delta}}_2]\hat u=r,
    ~{[{ {\hat \Delta}}_1~{ {\hat \Delta}}_2]}^H w=\hat s\right \}.
\end{eqnarray*}
Therefore, the proof is analogous to that of Theorem~\ref{thm:ustrJE} by using first Lemma~\ref{lem:JRBbstr_unstr}
and then Theorem~\ref{thm:spcl_matrix_map} for $1)$ and
Theorem~\ref{cor:real_spcl_cor} for $2)$.
\eproof

The following theorem presents the value of $\eta^{\mathcal S}(J,R,B,\lambda,x)$ and its
real counterpart if the original pencil is real. It also gives
$\eta^{\mathcal S}(J,R,B,\lambda):=\inf_{x\in\C^{2n+m}\setminus\{0\}}\eta^{\mathcal S}(J,R,B,\lambda,x)$.

\begin{theorem}\label{thm:strJRB}
Let $L(z)$ be a pencil as in~\eqref{eq:defL_cz}, and let
$\lambda \in i \R$ and $x \in \C^{2n+m}\setminus \{0\}$. Partition $x=[x_1^T~x_2^T~x_3^T]^T$
such that $x_1,\,x_2\in \C^n$ and $x_3 \in \C^{m}$ and define $ u=[x_2^T~~x_3^T]^T$, $w=x_1$,
\begin{eqnarray*}
r=(J-R+\lambda E)x_2 +Bx_3\quad \text{and}\quad s=\big[-((J+R+\lambda E)x_1)^T~~(B^Hx_1+Sx_3)^T\big]^T.
\end{eqnarray*}
Then the following statements hold.
\begin{enumerate}
\item $\eta^{\mathcal S}(J,R,B,\lambda,x)$ is finite if and only if $x_3=0$. In such a case the following holds.
\begin{eqnarray*}\label{eq:JRB_strberr}
 \eta^{{\mathcal S}}(J,R,B,\lambda,x) =
{\sqrt{{\|\widehat {\Delta}_1\|}_F^2+{\|\widehat{\Delta}_2\|}_F^2}},
\end{eqnarray*}
and
\begin{equation*}
 \eta^{\mathcal S}(J,R,B,\lambda)=\min\left\{\sigma_{\min}\left(\mat{cc}(J-R+\lambda E)&  B\rix^H\right),
 \sigma_{\min}(J-R+\lambda E)\right\},
\end{equation*}
where $\widehat{\Delta}_1\in\mathbb C^{n,n}$ and $\widehat{\Delta}_2\in\mathbb C^{n,m}$ are given by
\begin{equation*}\label{eq:str_JRBsepbrrobspecial}
[\widehat{\Delta}_1~\widehat{\Delta}_2]= \left\{\begin{array}{ll}
\frac{ru^H}{\|u\|^2} & \mbox{ if } x_1 = 0, \\[1ex]
\frac{ws^H}{\|w\|^2} & \mbox{ if } x_2 = 0,\\[1ex]
\frac{ru^H}{\|u\|^2} + \frac{ws^H}{\|w\|^2}\left(I_{n+m} - \frac{uu^H}{\|u\|^2}\right)
& \mbox{ otherwise. }  \end{array}\right.
\end{equation*}
\item If  $L(z)$ is real, and $\operatorname{rank}\left([x_1~\overline x_1]\right)=
\operatorname{rank}\left([x_2~\overline x_2]\right)=2$ then
$\eta^{{\mathcal S}_\R}(J,R,B,\lambda,x)$ is finite if and only if $x_3=0$ and $\lambda x_2^TEx_1=0$.
If the latter conditions are satisfied then
\begin{eqnarray*}\label{eq:JRBreal_strberr}
\eta^{{\mathcal S}_\R}(J,R,B,\lambda,x) =
{\sqrt{{\|\widetilde {\Delta}_1\|}_F^2+{\|\widetilde{\Delta}_2\|}_F^2}},
\end{eqnarray*}
where $\widetilde{\Delta}_1 \in \R^{n,n}$ and $\widetilde{\Delta}_2 \in \R^{n , m}$ are given by
\begin{equation*}\label{eq:str_JRBreal_sepbrrobspecial}
[\widetilde{\Delta}_1~\widetilde{\Delta}_2]= [r~ \overline r][u~ \overline u]^{\dag}+
\big([s~ \overline s][w~ \overline w]^{\dag}\big)^H
-\big([s~ \overline s][w~ \overline w]^{\dag}\big)^H\big([u~ \overline u][u~ \overline u]^{\dag}\big).
\end{equation*}
\end{enumerate}
\end{theorem}

\proof
The proof is similar to that of Theorem~\ref{thm:ustrJE} by using first Lemma~\ref{lem:JRBbstr_unstr}, and then
Theorem~\ref{thm:spcl_matrix_map} for $1)$ and
Theorem~\ref{cor:real_spcl_cor} for $2)$.
\eproof

\subsection{Perturbations to R, E, and B  or to J, E, and B}

This section is devoted to the block- and symmetry-structure-preserving eigenpair and eigenvalue backward
errors when only the blocks $R$, $E$ and $B$ of a pencil $L(z)$ as in~\eqref{eq:defL_cz}
are subject to perturbations. Let $\lambda \in \C$ and $x \in \C^{2n+m}\setminus \{0\}$, then
in view of Section~\ref{sec:section_def_epbr}, we have the definitions
\begin{eqnarray*}\label{eq:def_REB_eperr}
\eta^{\mathcal B}(R,E,B,\lambda,x)&=&\inf\Big\{{\|[\Delta_R~\Delta_E~\Delta_B]\|}_{F} \,\Big|~
\big((M-\Delta_M)+\lambda(N-\Delta_N)\big)x=0,\Delta_M + z \Delta_N \in \mathcal B\Big\},\\
\eta^{\mathcal S}(R,E,B,\lambda,x)&=&\inf\Big\{{\|[\Delta_R~\Delta_E~\Delta_B]\|}_F \,\Big |~
\big((M-\Delta_M)+\lambda(N-\Delta_N)\big)x=0,\Delta_M + z \Delta_N \in \mathcal S\Big\},
\end{eqnarray*}
respectively, where $\mathcal B$ is the set of all pencils
of the form $\Delta L(z)=\Delta_M + z \Delta_N$ with
\[
\Delta_M= \mat{ccc} 0 &-\Delta_R & \Delta_B\\ {-\Delta_R }^H &0&0 \\\Delta_B^H&0&0\rix
\quad \text{and}\quad \Delta_N=\mat{ccc} 0 &\Delta_E & 0\\ -\Delta_E ^H &0&0 \\0&0&0\rix
\]
and $\Delta R,\Delta E\in\mathbb C^{n,n}$, $\Delta B\in\mathbb C^{n,m}$, and $\mathcal S$ is the
corresponding set of all such pencils that in addition satisfy $\Delta R,\Delta E\in\text{Herm}(n)$.

\begin{remark}\rm
Let $L(z)$ be a pencil as in~\eqref{eq:defL_cz}, and let
$\lambda \in\mathbb C$ and $x=[x_1^T~x_2^T~x_3^T]^T$ be such that $x_1,x_2 \in \C^n$ and $x_3 \in \C^m$.
Then for any $\Delta_R,\, \Delta_E\in \C^{n,n}$ and $\Delta_B \in \C^{n , m}$,
and corresponding $\Delta L(z)=\Delta_M +z \Delta_N \in \mathcal B$, we have
$(L-\Delta L)(\lambda)x=0$ if and only if
\begin{eqnarray*}
(-\Delta_R+\lambda \Delta_E) x_2 + \Delta_B x_3 &=&(J-R+\lambda E)x_2 +B x_3, \\
(-\Delta_R+\lambda \Delta_E)^Hx_1 &=&{(-J-R-\lambda E)x_1},\\
\Delta_B^Hx_1&=& B^Hx_1+Sx_3,
\end{eqnarray*}
which, in turn, is equivalent to
\begin{eqnarray}
 \mat{cc} -\Delta_R+\lambda \Delta_E &\Delta_B \rix \underbrace{\mat{c} x_2\\x_3 \rix}_{=u} &=&
 \underbrace{(J-R+\lambda E)x_2 +B x_3}_{=r}, \label{REBequi_1}\\
 \mat{cc} -\Delta_R+\lambda \Delta_E &\Delta_B \rix^H \underbrace{x_1}_{=w} &=&
 \underbrace{\mat{c} -(J+R+\lambda E)x_1\\ B^Hx_1 + Sx_3 \rix}_{=s}. \label{REBequi_2}
\end{eqnarray}
\end{remark}

\begin{lemma}\label{lem:REBbstr_unstr}
Let $L(z)$ be a pencil as in~\eqref{eq:defL_cz}, and let
$\lambda \in i \R$ and $x \in \C^{2n+m}\setminus \{0\}$. Partition $x=[x_1^T~x_2^T~x_3^T]^T$
such that $x_1,\,x_2 \in \C^n$ and $x_3 \in \C^{m}$ and let $u\,,w\,,r$ and $s$ be defined as
in~\eqref{REBequi_1} and~\eqref{REBequi_2}. Then the following statements are equivalent.
\begin{enumerate}
 \item  There exist $\Delta_R,\,\Delta_E\in \C^{n,n}$ and $\Delta_B \in \C^{n , m}$
  satisfying~\eqref{REBequi_1} and~\eqref{REBequi_2}.
 \item There exist ${\Delta} \in \C^{n , n+m}$ such that $\Delta u=r$ and ${\Delta}^H w=s$.
  \item  There exist $\Delta_R,\,\Delta_E \in {\rm Herm}(n)$ and
 $\Delta_B \in \C^{n , m}$ satisfying~\eqref{REBequi_1} and~\eqref{REBequi_2}.
 \item  $x$ satisfies $x_3=0$.
\end{enumerate}
Moreover,
\begin{eqnarray*}\label{REBequi_3}
 && \inf\left\{ {\|\Delta_R ~ \Delta_E~\Delta_B\|}_F^2 \,\Big|~ \Delta_R, \, \Delta_E\in \C^{n,n},\,
  \Delta_B \in \C^{n , m} ~{\rm satisfy}~\eqref{REBequi_1}\,{\rm and}~\eqref{REBequi_2} \right \} \\
&&\ =\inf \left\{\left.\frac{{\|{\Delta}_1\|}_F^2}{1+|\lambda|^2}+
  {\|{\Delta}_2\|}_F^2\,\right |~ {\Delta}_1 \in \C^{n,n},\,{\Delta}_2 \in \C^{n , m},\,
  [{\Delta}_1~{\Delta}_2] u=r,~[{\Delta}_1~{\Delta}_2] ^H w=s\right \},
 \end{eqnarray*}
and
\begin{eqnarray*}\label{REBequi_4}
&&  \inf\left\{ {\|\Delta_R ~ \Delta_E~\Delta_B\|}_F^2 \,\Big|~ \Delta_R, \, \Delta_E \in {\rm Herm}(n),\,
  \Delta_B \in \C^{n , m} ~{\rm satisfy}~\eqref{REBequi_1}\,{\rm and}~\eqref{REBequi_2} \right \} \\
&&\ =  \inf \Bigg\{{\left \|\frac{{\Delta}_1 + {\Delta}_1^H}{2} \right\|}_F^2 +
  \frac{1}{|\lambda|^2}{\left \|\frac{{\Delta}_{1} - {\Delta}_1^H}{2} \right\|}_F^2+
  {\|{\Delta}_2\|}_F^2~\Bigg |~ {\Delta}_1 \in \C^{n,n},\,  {\Delta}_2 \in \C^{n,m},\,\\
&& \phantom{\ =\inf \Bigg\{{\left \|\frac{{\Delta}_1 + {\Delta}_1^H}{2} \right\|}_F^2 +
  \frac{1}{|\lambda|^2}{\left \|\frac{{\Delta}_{1} - {\Delta}_1^H}{2} \right\|}_F^2+
  {\|{\Delta}_2\|}_F^2~\Bigg |~} \,[{\Delta}_1~{\Delta}_2] u=r,~[{\Delta}_1~{\Delta}_2] ^H w=s\Bigg \}.
\end{eqnarray*}
\end{lemma}
\proof
Again, by using the fact that any $\Delta \in \C^{n , n+m}$ can be written as
$\Delta=[{\Delta}_1\,~{\Delta}_2]$ where ${\Delta}_1\in \C^{n,n}$ and ${\Delta}_2\in \C^{n , m}$
such that ${\|\Delta\|}_F={\|[{\Delta}_1\,~{\Delta}_2]\|}_F
=\sqrt{{{\|{\Delta}_1\|}_F^2+{\|{\Delta}_2\|}}_F^2}$, the
proof is obtained by arguments similar to those in the proof Lemma~\ref{thm;equivalentcond}
and Lemma~\ref{lem1:str_RE}.
\eproof

\begin{theorem} \label{thm:ustrREB}
Let $L(z)$ be a pencil as in~\eqref{eq:defL_cz}, and let
$\lambda \in i \R$ and $x \in \C^{2n+m}\setminus \{0\}$. Partition $x=[x_1^T~x_2^T~x_3^T]$
so that $x_1,\,x_2\in \C^n$, and $x_3 \in \C^{m}$, and define $w=x_1$,
$u=[x_2^T~x_3^T]^T$, $\hat u=[{(1+|\lambda|^2)^{1/2}}x_2^T~~x_3^T]^T$,
$r=(J-R+\lambda E)x_2 +Bx_3$,  $
 s=[-((J+R+\lambda E)x_1)^T~(B^Hx_1+Sx_3)^T]^T$,
and $\widehat s=\big[-(1+|\lambda|^2)^{-1/2}((J+R+\lambda E)x_1)^T~~(B^Hx_1+Sx_3)^T\big]^T$.
Then
$\eta^{\mathcal B}(R,E,B,\lambda,x)$ and $\eta^{\mathcal S}(R,E,B,\lambda,x)$
are finite if and only if $x_3=0$. Furthermore, the following statements hold.
\begin{enumerate}
\item If $x_3=0$, then
\begin{eqnarray}\label{eq:REB_unstrberr}
 \eta^{\mathcal B}(R,E,B,\lambda,x) ={\sqrt{{\|\widehat {\Delta}_1\|}_F^2+{\|\widehat{\Delta}_2\|}_F^2}},
\end{eqnarray}
and
\begin{equation}\label{eq:REB24aug}
 \eta^{\mathcal B}(R,E,B,\lambda)=\min\left\{\sigma_{\min}\left(\mat{cc}\frac{(J-R+\lambda E)}{\sqrt{1+|\lambda|^2}}&  B\rix^H\right),
 \frac{\sigma_{\min}(J-R+\lambda E)}{\sqrt{1+|\lambda|^2}}\right\},
\end{equation}
where $\widehat{\Delta}_1$ and $\widehat{\Delta}_2$ are given by
\begin{equation*}\label{eq:unstrREBsepbrrobspecial}
[\widehat{\Delta}_1~\widehat{\Delta}_2]= \left\{\begin{array}{ll}
\frac{r\hat u^H}{\|\hat u\|^2} & \mbox{ if } x_1 = 0, \\[1ex]
\frac{w \hat s^H}{\|w\|^2} & \mbox{ if } x_2 = 0,\\[1ex]
\frac{r\hat u^H}{\|\hat u\|^2} + \frac{w\hat s^H}{\|w\|^2}\left(I_{n+m} - \frac{\hat u \hat u^H}{\|\hat u\|^2}\right)
& \mbox{ otherwise. }  \end{array}\right.
\end{equation*}
\item If $x_3=0$, then
{\small{\begin{eqnarray}\label{eq:REB_strberrbound1}
{\sqrt{{{\|\widetilde {\Delta}_1\|}_F^2}+{\|\widetilde{\Delta}_2\|}_F^2}}
\leq  \eta^{{\mathcal S}}(R,E,B,\lambda,x) \leq
\sqrt{ {\left\|\frac{{\widetilde {\Delta}}_1 + {\widetilde {\Delta}}_1^H}{2} \right\|}_F^2 +
  \frac{1}{|\lambda|^2}{\left \|\frac{{\widetilde {\Delta}}_{1} - {\widetilde {\Delta}}_1^H}{2} \right\|}_F^2+
  {\|{\widetilde {\Delta}}_2\|}_F^2}, \qquad
\end{eqnarray}}}
when $|\lambda|\leq 1$, and
{\small{\begin{eqnarray}\label{eq:REB_strberrbound2}
{\sqrt{\frac{{\|\widetilde {\Delta}_1\|}_F^2}{|\lambda|^2}+{\|\widetilde{\Delta}_2\|}_F^2}}
\leq \eta^{{\mathcal S}}(R,E,B,\lambda,x) \leq
\sqrt{{\left\|\frac{{\widetilde {\Delta}}_1 + {\widetilde {\Delta}}_1^H}{2} \right\|}_F^2 +
  \frac{1}{|\lambda|^2}{\left \|\frac{{\widetilde {\Delta}}_{1} - {\widetilde {\Delta}}_1^H}{2} \right\|}_F^2+
  {\|{\widetilde {\Delta}}_2\|}_F^2},\qquad
\end{eqnarray}}}
when $|\lambda|\geq 1$, where $\widetilde{\Delta}_1$ and $\widetilde{\Delta}_2$ are given by
\begin{equation*}\label{eq:strREBsepbrrobspecial}
[\widetilde{\Delta}_1~\widetilde{\Delta}_2]= \left\{\begin{array}{ll}
\frac{ru^H}{\|u\|^2} & \mbox{ if } x_1 = 0, \\[1ex]
\frac{ws^H}{\|w\|^2} & \mbox{ if } x_2 = 0,\\[1ex]
\frac{ru^H}{\|u\|^2} + \frac{ws^H}{\|w\|^2}\left(I_{n+m} - \frac{uu^H}{\|u\|^2}\right)
& \mbox{ otherwise. }  \end{array}\right.
\end{equation*}
\end{enumerate}
\end{theorem}
\proof
In view of~\eqref{REBequi_1} and~\eqref{REBequi_2}, we have
\begin{eqnarray*}\label{eq:str_REEupdated1}
&&\big(\eta^{\mathcal B}(R,E,B,\lambda,x)\big)^2\\
&&\qquad=\inf\Big\{{\left\|[\Delta_R~\Delta_E~\Delta_B]\right\|}_{F}^2 \,\Big |
\Delta_R,\,\Delta_E\in \C^{n,n},\,\Delta_B \in \C^{n , m}\;
 {\rm satisfy}\,\eqref{REBequi_1}~{\rm and}~\eqref{REBequi_2}\Big\} \\
&&\qquad=\inf \bigg\{\frac{{\|{\Delta}_1\|}_F^2}{1+|\lambda|^2}+
  {\|{\Delta}_2\|}_F^2\,\bigg |~ {\Delta}_1\in \C^{n,n},\,{\Delta}_2 \in \C^{n , m},\,
    [{\Delta}_1~{\Delta}_2] u=r,~[{\Delta}_1~{\Delta}_2] ^H w=s\bigg \},
\end{eqnarray*}
where the last equality follows from  Lemma~\ref{lem:REBbstr_unstr}. Observe that if we set
$\widehat {\Delta}_2={\Delta}_2$ and
$\widehat {\Delta}_1=\frac{{\Delta}_1}{\sqrt{1+|\lambda|^2}}$, then we obtain
\begin{eqnarray*}
&&\big(\eta^{\mathcal B}(R,E,B,\lambda,x)\big)^2\\
&&\qquad=\inf \Big\{{\|\widehat{\Delta}_1\|}_F^2+
  {\|\widehat{\Delta}_2\|}_F^2~\Big |~ \widehat{\Delta}_1\in \C^{n,n},\,\widehat{\Delta}_2 \in \C^{n , m},\,
  [\widehat{\Delta}_1~\widehat{\Delta}_2] \widehat u=r,
  ~[\widehat{\Delta}_1~\widehat{\Delta}_2] ^H w=\widehat s \Big \}.
\end{eqnarray*}
Thus~\eqref{eq:REB_unstrberr} follows from Theorem~\ref{thm:spcl_matrix_map}, and arguments similar to
those in the proof of Theorem~\ref{thm:ustrJE} give~\eqref{eq:REB24aug}.

Similarly, by using Lemma~\ref{lem:REBbstr_unstr} in the definition of $\eta^{\mathcal S}(R,E,B,\lambda,x)$
we can write
\begin{eqnarray*}
\big(\eta^{\mathcal S}(R,E,B,\lambda,x)\big)^2 =  \inf \Bigg\{{\left \|\frac{{\Delta}_1 + {\Delta}_1^H}{2}
\right\|}_F^2 +  \frac{1}{|\lambda|^2}{\left\|\frac{{\Delta}_{1} - {\Delta}_1^H}{2} \right\|}_F^2+
  {\|{\Delta}_2\|}_F^2~\Bigg |~ {\Delta}_1,\,{\Delta}_2 \in \C^{n , n},\\
  \,[{\Delta}_1~{\Delta}_2] u=r,~[{\Delta}_1~{\Delta}_2] ^H w=s\Bigg \}. \nonumber
\end{eqnarray*}
For any ${\Delta}_1 \in \C^{n , n}$ we have
${\|{\Delta}_1\|}_F^2={\left\|\frac{{\Delta}_1 +{\Delta}_1^H}{2}\right\|}_F^2+
{\left\|\frac{{\Delta}_1 -{\Delta}_1^H}{2}\right\|}_F^2$. This implies
\begin{eqnarray}\label{eq:REB_auxi_strberrbound1}
{\|{\Delta}_1\|}_F^2+{\|{\Delta}_2\|}_F^2~\leq ~ {\left\|\frac{{\Delta}_1 +{\Delta}_1^H}{2}\right\|}_F^2+
\frac{1}{|\lambda|^2}{\left\|\frac{{\Delta}_1 -{\Delta}_1^H}{2}\right\|}_F^2+
{\|{\Delta}_2\|}_F^2~ \quad \quad {\rm if }~|\lambda|\leq 1
\end{eqnarray}
and
\begin{eqnarray}\label{eq:REB_auxi_strberrbound2}
\frac{{\| {\Delta}_1\|}_F^2}{|\lambda|^2}+{\|{\Delta}_2\|}_F^2~\leq ~ {\left\|\frac{{\Delta}_1 +{\Delta}_1^H}{2}\right\|}_F^2+
\frac{1}{|\lambda|^2}{\left\|\frac{{\Delta}_1 -{\Delta}_1^H}{2}\right\|}_F^2
+{\|{\Delta}_2\|}_F^2~
\quad \quad {\rm if }~|\lambda|\geq 1
\end{eqnarray}
for all $\Delta_1\in \C^{n,n}$ and $\Delta_2 \in \C^{n , m}$. Taking the infimum over all
${\Delta}_1 \in \C^{n,n}$, ${\Delta}_2 \in \C^{n , m}$ satisfying $[{\Delta}_1~{\Delta}_2] u=r$ and
$[{\Delta}_1~{\Delta}_2]^Hw=s$ in~\eqref{eq:REB_auxi_strberrbound1} and~\eqref{eq:REB_auxi_strberrbound2}
followed by applying Theorem~\ref{thm:spcl_matrix_map} yields~\eqref{eq:REB_strberrbound1}
and~\eqref{eq:REB_strberrbound2}.
\eproof

\begin{remark}\label{thm:ustrJEB}
{\rm
 We mention that a result similar to Theorem \ref{thm:ustrREB}  can also be obtained for
 the block-structure-preserving eigenpair and eigenvalue backward errors $\eta^{\mathcal B}(J,E,B,\lambda,x)$ and
 $\eta^{\mathcal B}(J,E,B,\lambda)$, respectively, when perturbations are restricted to
 affect only the blocks $J$, $E$ and $B$ of a pencil $L(z)$ as in~\eqref{eq:defL_cz}. In fact,
 for $\lambda\in i\mathbb R$ and $x\in\mathbb C^{2n+m}$, using arguments analogous to those
 in this section we obtain that
 \[
  \eta^{\mathcal B}(J,E,B,\lambda,x)=\eta^{\mathcal B}(R,E,B,\lambda,x) \quad \text{and}\quad
  \eta^{\mathcal B}(J,E,B,\lambda)=\eta^{\mathcal B}(R,E,B,\lambda).
 \]
 }
\end{remark}

\subsection{Perturbation to J, R, and E }

Let $L(z)$ be a pencil as in~\eqref{eq:defL_cz}, and let $\lambda \in \C$ and $x \in \C^{2n+m}\setminus \{0\}$.
In this section, we allow perturbations in the blocks $J$, $R$ and $E$ of $L(z)$. The
block- and symmetry-structure-preserving eigenpair backward errors $\eta^{\mathcal B}(J,R,E,\lambda,x)$ and
$\eta^{\mathcal S}(J,R,E,\lambda,x)$ are defined by
\begin{eqnarray*}\label{eq:def_JRE_eperr}
\eta^{\mathcal B}(J,R,E,\lambda,x)&=&\inf\Big\{\big\|[\Delta_J~\Delta_R~\Delta_E]\big\|_{F} \Big |
\big((M-\Delta_M)+\lambda(N-\Delta_N)\big)x=0,\,\Delta_M + z \Delta_N \in \mathcal B\Big\},\\
\eta^{\mathcal S}(J,R,E,\lambda,x)&=&\inf\Big\{\big\|[\Delta_J~\Delta_R~\Delta_E]\big\|_{F} \Big |
\big((M-\Delta_M)+\lambda(N-\Delta_N)\big)x=0,\,\Delta_M + z \Delta_N \in \mathcal S\Big\},
\end{eqnarray*}
respectively, where $\mathcal B$ is the set of all pencils of the form $\Delta L(z)=
\Delta_M + z \Delta_N$ with
\[
\Delta_M= \mat{ccc} 0 &\Delta_J-\Delta_R & 0 \\ {(\Delta_J-\Delta_R )}^H &0&0 \\0&0&0\rix
\quad \text{and}\quad
\Delta_N=\mat{ccc} 0 &\Delta_E & 0\\ -\Delta_E^H &0&0 \\0&0&0\rix
\]
and $\mathcal S$ is the corresponding set of pencils from $\mathcal B$ that satisfy in addition
$\Delta_J\in\text{SHerm}(n)$ and $\Delta_R,\Delta_E\in\text{Herm}(n)$.

\begin{remark}\rm
If $\lambda \in \mathbb C$ and $x=[x_1^T~x_2^T~x_3^T]^T$ are such that $x_1,x_2 \in \C^n$ and $x_3 \in \C^m$,
then for any $\Delta_J,\, \Delta_R,\,\Delta_E \in \C^{n , n}$ and corresponding
$\Delta L(z)=\Delta_M +z \Delta_N \in \mathcal B$, we have $(L-\Delta L)(\lambda)x=0$ if and only if
\begin{eqnarray}
(\Delta_J-\Delta_R+\lambda \Delta_E) \underbrace{x_2}_{=u}&=&\underbrace{(J-R+\lambda E)x_2 +B x_3}_{=r}, \label{JREequi_1}\\
(\Delta_J-\Delta_R+\lambda \Delta_E)^H\underbrace{x_1}_{=w} &=&\underbrace{{(-J-R-\lambda E)x_1}}_{=s},\label{JREequi_2}\\
0&=& B^Hx_1+Sx_3. \label{JREequi_3}
\end{eqnarray}
\end{remark}

\begin{lemma}\label{lem:JREequivalentcond}
Let $L(z)$ be a pencil as in~\eqref{eq:defL_cz}, and let $\lambda \in i \R$ and $x \in \C^{2n+m}\setminus \{0\}$.
Partition $x=[x_1^T~x_2^T~x_3^T]^T$ such that $x_1,\,x_2\in \C^n$ and $x_3 \in \C^{m}$ and let
$u\,,w\,,r$ and $s$ be defined as in~\eqref{JREequi_1} and~\eqref{JREequi_2}.
 Then the following statements are equivalent.
\begin{enumerate}
 \item There exist $\Delta_J,\Delta_R,\,\Delta_E \in \C^{n , n}$ satisfying~\eqref{JREequi_1} and~\eqref{JREequi_2}.
 \item There exist $\Delta \in \C^{n , n}$ such that $\Delta u=r$ and ${\Delta}^H w=s$.
 \item There exist $\Delta_J \in {\rm SHerm}(n)$ and $\Delta_R,\,\Delta_E \in {\rm Herm}(n)$
 satisfying~\eqref{JREequi_1} and~\eqref{JREequi_2}.
 \item  $x$ satisfies $x_3^HB^Hx_1=0$.
\end{enumerate}
Moreover, we have
\begin{eqnarray}\label{JREequi_4}
  \inf\Big\{ {\big\|\big[\Delta_J~\Delta_R~\Delta_E\big]\big\|}_F^2  ~\Big|~ \Delta_J,\,\Delta_R,\,\Delta_E
  \in \C^{n , n} ~{\rm satisfy}~\eqref{JREequi_1}~{\rm and}~\eqref{JREequi_2}  \Big \} \nonumber \\
  =\inf \left\{\left.\frac{{\|\Delta\|}_F^2}{2+|\lambda|^2}~\right|~ \Delta \in \C^{n , n},\,
  \Delta u=r,~\Delta ^H w=s\right\},
\end{eqnarray}
and
\begin{eqnarray}\label{JREequi_5}
  \inf\left\{\left. {\big\|\big[\Delta_J~\Delta_R~\Delta_E\big]\big\|}_F^2  \,\right|\, \Delta_J \in {\rm SHerm}(n),\,
  \Delta_E,\, \Delta_R\in {\rm Herm}(n) ~{\rm satisfying}~\eqref{JREequi_1}\,{\rm and}~\eqref{JREequi_2}\right \} \nonumber \\
 = \inf \left\{\left.{\left \|\frac{\Delta + \Delta^H}{2} \right\|}_F^2 +
  \frac{1}{1+|\lambda|^2}{\left \|\frac{\Delta - \Delta^H}{2} \right\|}_F^2\,\right|~ \Delta \in \C^{n , n},\,
  \Delta u=r,~\Delta ^H w=s\right \}.\qquad
\end{eqnarray}
\end{lemma}

\proof
``$1) \Rightarrow 2)$'': Let $\Delta_J$, $\Delta_R$, $\Delta_E\in \C^{n , n}$ be such that they
satisfy~\eqref{JREequi_1} and~\eqref{JREequi_2}. By setting $\Delta =\Delta_J-\Delta_R + \lambda \Delta_E$
we get $\Delta u=r$ and ${\Delta}^Hw=s$. Also, we obtain
\begin{equation}\label{JREequi_6}
 {\|\Delta\|}_F^2 \leq \big({\|\Delta_J\|}_F+{\|\Delta_R\|}_F+|\lambda|{\|\Delta_E\|}_F \big)
  \leq (2+|\lambda|^2)\big({\| \Delta_J\|}_F^2+{\|\Delta_R\|}_F^2+{\| \Delta_E\|}_F^2\big),
\end{equation}
where the latter inequality follows from the Cauchy-Schwarz inequality (in $\mathbb R^3$).
Then ``$\geq$'' in~\eqref{JREequi_4} can be shown similarly as ``$1)\Rightarrow2)$''
in the proof of Lemma~\ref{thm;equivalentcond}.

\noindent
``$2) \Rightarrow 1)$'': Conversely, let $\Delta \in \C^{n , n}$ such that $\Delta u=r$ and ${\Delta}^Hw=s$.
Define
\[
\Delta_J=\frac{\Delta}{2+|\lambda|^2},\quad \Delta_R=-\frac{\Delta}{2+|\lambda|^2},\quad\mbox{and}\quad
\Delta_E= \frac{\bar \lambda \Delta}{2+|\lambda|^2}.
\]
Then $\Delta_J,\,\Delta_R$ and $\Delta_E$ satisfy $\Delta_J -\Delta_R+\lambda \Delta_E =\Delta$ and
hence~\eqref{JREequi_1} and~\eqref{JREequi_2}. Furthermore, we have
\[
 {\| \Delta_J\|}_F^2+{\|\Delta_R\|}_F^2+{\| \Delta_E\|}_F^2= \frac{{\|\Delta\|}_F^2}{2+|\lambda|^2}.
\]
Thus, we get ``$\leq$'' in~\eqref{JREequi_4} by following arguments similar to those of ``$2) \Rightarrow 1)$''
in the proof of Lemma~\ref{thm;equivalentcond}.

\noindent
``$2)\Rightarrow 3)$'': To show this, let $\Delta \in \C^{n , n}$ be such that $\Delta u=r$ and ${\Delta}^Hw=s$.
Then by setting
\[
\Delta_R=-\frac{\Delta+\Delta^H}{2},\quad \Delta_J=\frac{\Delta-\Delta^H}{2(1+|\lambda|^2)},\quad\mbox{and}\quad
\Delta_E= \frac{\bar \lambda( \Delta-\Delta^H)}{2(1+|\lambda|^2)},
\]
we have $\Delta_J \in {\rm SHerm}(n)$ and
$\Delta_R,\,\Delta_E \in {\rm Herm}(n)$ (using $\lambda \in i\R$),  and furthermore we obtain
\[
\Delta_J -\Delta_R +\lambda \Delta_E = \frac{\Delta-{\Delta}^H}{2(1+|\lambda|^2)}
+ \frac{\Delta+{\Delta}^H}{2} + \frac{{|\lambda|}^2( \Delta-{\Delta}^H)}{2(1+|\lambda|^2)}
= \frac{\Delta+{\Delta}^H}{2} +\frac{\Delta-{\Delta}^H}{2} = \Delta.
\]
Thus, $\Delta_J$, $\Delta_R$, and $\Delta_E$ satisfy~\eqref{JREequi_1} and~\eqref{JREequi_2}, and also
\[
 {\| \Delta_J\|}_F^2+{\|\Delta_R\|}_F^2+{\| \Delta_E\|}_F^2=
 {\left \|\frac{\Delta + \Delta^H}{2} \right\|}_F^2 +
  \frac{1}{1+|\lambda|^2}{\left \|\frac{\Delta - \Delta^H}{2} \right\|}_F^2.
\]
Now ``$\leq$'' in~\eqref{JREequi_5} can be shown by arguments similar to those of ``$2) \Rightarrow 1)$''
in the proof of Lemma~\ref{lem1:str_RE}.

\noindent
``$3) \Rightarrow 2)$'': Suppose that $\Delta_J \in {\rm SHerm}(n)$ and $\Delta_R,\,\Delta_E \in {\rm Herm}(n)$
satisfy~\eqref{JREequi_1} and~\eqref{JREequi_2}. Define $\Delta = \Delta_J -\Delta_R+\lambda \Delta_E$,
 then $\Delta u=r$ and ${\Delta }^Hw=s$. Note that $\Delta_J +\lambda \Delta_E$ is skew-Hermitian
 since $\lambda \in i\R$, and therefore  $\Delta_J +\lambda \Delta_E$ and $-\Delta_R$ are respectively
 the unique skew-Hermitian and Hermitian parts of $\Delta$, i.e.,
 \[
 \Delta_R= - \frac{\Delta +\Delta^H}{2}\quad\mbox{and}\quad
 \Delta_J +\lambda \Delta_E=\frac{\Delta -\Delta^H}{2}.
 \]
This implies
\[
\left \|\frac{\Delta -{\Delta}^H}{2} \right \|_F=
\|\Delta_J+ \lambda \Delta_E\|_F \leq \|\Delta_J\|_F +|\lambda|\cdot\|\Delta_E\|_F
 \]
 and
\[
{\frac{1}{1+|\lambda|^2}}{\left \|\frac{\Delta -{\Delta}^H}{2} \right \|}_F^2
\leq {\|\Delta_J\|}_F^2+{\|\Delta_E\|}_F^2,
\]
where the last inequality is obtained with the help of the Cauchy-Schwarz inequality (in $\mathbb R^2$).
Furthermore, we have
\begin{equation}\label{JREequi_7}
{\frac{1}{1+|\lambda|^2}}{\left \|\frac{\Delta -{\Delta}^H}{2} \right \|}_F^2
+{\left \|\frac{\Delta +{\Delta}^H}{2} \right \|}_F^2
\leq {\|\Delta_J\|}_F^2+{\|\Delta_E\|}_F^2 +{\|\Delta_R\|}_F^2.
\end{equation}
Thus, arguments similar to those in ``$1) \Rightarrow 2)$'' in the proof of Lemma~\ref{thm;equivalentcond}
give ``$\geq$'' in~\eqref{JREequi_5}.

\noindent
``$2)\Leftrightarrow 4)$'': This follows immediately from Theorem~\ref{thm:spcl_matrix_map}.
\eproof

\begin{theorem} \label{thm:ustrJRE}
 Let $L(z)$ be a pencil as in~\eqref{eq:defL_cz}, and let
 $\lambda \in i \R$ and $x \in \C^{2n+m}\setminus \{0\}$. Partition $x=[x_1^T~x_2^T~x_3^T]^T$
 such that $x_1,\,x_2 \in \C^n$ and $x_3 \in \C^{m}$ and define
$r=(J-R+\lambda E)x_2 +Bx_3$ and $s=-(J+R+\lambda E)x_1$.
Then
$\eta^{\mathcal B}(J,R,E,\lambda,x)$ and $\eta^{\mathcal S}(J,R,E,\lambda,x)$
are finite if and only if $x_3=0$ and $B^Hx_1=0$. If latter conditions are satisfied then
\begin{eqnarray}\label{eq:JRE_unstrberr}
 \eta^{\mathcal B}(J,R,E,\lambda,x) = \frac{\|\widehat{\Delta}\|_F}{\sqrt{2+|\lambda|^2}},\quad
 \eta^{\mathcal B}(J,R,E,\lambda)=\frac{\sigma_{\min}(J-R+\lambda E)}{\sqrt{2+|\lambda|^2}}
\end{eqnarray}
and
\begin{eqnarray}\label{eq:JRE_strberr}
\frac{\|\widehat{\Delta}\|_F}{\sqrt{1+|\lambda|^2}}\leq \eta^{\mathcal S}(J,R,E,\lambda,x) \leq
\sqrt{{\left \|\frac{\widehat{\Delta} + \widehat{\Delta}^H}{2} \right\|}_F^2 +
  \frac{1}{1+|\lambda|^2}{\left \|\frac{\widehat{\Delta} - \widehat{\Delta}^H}{2} \right\|}_F^2},
\end{eqnarray}
where $\widehat \Delta$ is given by
\begin{equation*}\label{eq:JREsepbrrobspecial}
\widehat \Delta= \left\{\begin{array}{ll}
\frac{ru^H}{\|u\|^2} & \mbox{ if } x_1 = 0, \\[1ex]
\frac{ws^H}{\|w\|^2} & \mbox{ if } x_2 = 0,\\[1ex]
\frac{ru^H}{\|u\|^2} + \frac{ws^H}{\|w\|^2}\left(I_{n} - \frac{uu^H}{\|u\|^2}\right)
& \mbox{ otherwise. }  \end{array}\right.
\end{equation*}
\end{theorem}

\proof In view of Lemma~\ref{lem:JREequivalentcond} and Theorem~\ref{thm:spcl_matrix_map}, the
proofs of~\eqref{eq:JRE_unstrberr} and~\eqref{eq:JRE_strberr} are based on similar arguments as those in the
proofs of Theorem~\ref{thm:ustrJE} and Theorem~\ref{thm:strRE}, respectively.
\eproof

\section{Perturbations in $J$, $R$, $E$, and $B$}

Finally, in this section we allow all four blocks $J$, $R$, $E$, and $B$ of a pencil $L(z)$
as in~\eqref{eq:defL_cz} to be perturbed. Let $\lambda \in \C$ and $x \in \C^{2n+m} \setminus \{0\}$, then
by the terminology of Section~\ref{sec:section_def_epbr} the block- and symmetry-structure-preserving
eigenpair backward errors $\eta^{\mathcal B}(J,R,E,B,\lambda,x)$ and
$\eta^{\mathcal S}(J,R,E,B,\lambda,x)$ are respectively defined by
\begin{eqnarray*}\label{eq:def_JREB_eperr}
&&\eta^{\mathcal B}(J,R,E,B,\lambda,x)\\
&&\qquad=\inf\Big\{{\big\|[\Delta_J~\Delta_R~\Delta_E~\Delta_B]\big\|}_{F}\, \Big |~
\big((M-\Delta_M)+\lambda(N-\Delta_N)\big)x=0,\,\Delta_M + z \Delta_N \in \mathcal B\Big\},\\
&&\eta^{\mathcal S}(J,R,E,B,\lambda,x)\\
&&\qquad=\inf\Big\{{\big\|[\Delta_J~ \,\Delta_R~\Delta_E~\Delta_B]\big\|}_{F} \,\Big |~
\big((M-\Delta_M)+\lambda(N-\Delta_N)\big)x=0,\,\Delta_M + z \Delta_N \in \mathcal S\Big\},
\end{eqnarray*}
where $\mathcal B$ denotes the set of all pencils of the
form $\Delta L(z)=\Delta_M + z \Delta_N$ with
\[
\Delta_M= \mat{ccc} 0 &\Delta_J-\Delta_R & \Delta_B \\ {(\Delta_J-\Delta_R )}^H &0&0 \\\Delta_B^H&0&0\rix
\quad \text{and}\quad
\Delta_N=\mat{ccc} 0 &\Delta_E & 0\\ -\Delta_E^H &0&0 \\0&0&0\rix
\]
and $\mathcal S$ is the set of corresponding pencils where in addition we have that $\Delta_J\in\text{SHerm}(n)$
and $\Delta_R,\Delta_E\in\text{Herm}(n)$.

\begin{remark}\rm
If $\lambda \in \mathbb C$ and $x=[x_1^T~x_2^T~x_3^T]^T $ are such that $x_1,x_2 \in \C^n$ and $x_3 \in \C^m$,
then for any $\Delta_B \in \C^{n,m}$, $\Delta_J,\,\Delta_R,\, \Delta_E \in \C^{n,n}$,
 and corresponding $\Delta L(z)=\Delta_M +z \Delta_N \in \mathcal B$, we have
$(L-\Delta L)(\lambda)(x)=0$ if and only if
\begin{eqnarray*}
(\Delta_J-\Delta_R+\lambda \Delta_E) x_2 + \Delta_B x_3 &=&(J-R+\lambda E)x_2 +B x_3, \\
(\Delta_J-\Delta_R+\lambda \Delta_E)^Hx_1 &=&{(-J-R-\lambda E)x_1},\\
\Delta_B^Hx_1&=& B^Hx_1+Sx_3,
\end{eqnarray*}
which in turn is equivalent to
\begin{eqnarray}
 \mat{cc} \Delta_J-\Delta_R+\lambda \Delta_E &\Delta_B \rix \underbrace{ \mat{c} x_2\\x_3 \rix}_{=u}
  &=& \underbrace{(J-R+\lambda E)x_2 +B x_3}_{=r}, \label{JREBequi_1} \\
 \mat{cc} \Delta_J-\Delta_R+\lambda \Delta_E &\Delta_B \rix^H \underbrace{x_1}_{=w} &=&\underbrace{\mat{c} -(J+R+\lambda E)x_1\\ B^Hx_1 + Sx_3 \rix}_{=s}. \label{JREBequi_2}
\end{eqnarray}
\end{remark}

\begin{lemma}\label{lem:JREBequivalentcond}
Let $L(z)$ be a pencil as in~\eqref{eq:defL_cz}, and let
$\lambda \in i \R$ and $x \in \C^{2n+m}\setminus \{0\}$. Partition $x=[x_1^T~x_2^T~x_3^T]^T$
 such that $x_1,\,x_2\in\C^n$ and $x_3 \in \C^{m}$ and let $u\,,w\,,r$ and $s$ be defined as in~\eqref{JREBequi_1}
 and~\eqref{JREBequi_2}. Then the following statements are equivalent.
\begin{enumerate}
 \item There exist $\Delta_J,\,\Delta_R,\,\Delta_E\in \C^{n,n}$ and $\Delta_B \in \C^{n , m}$
 satisfying~\eqref{JREBequi_1} and~\eqref{JREBequi_2}.
 \item There exist ${\Delta} \in \C^{n , n+m}$ such that ${\Delta}u=r$ and ${\Delta}^H w=s$.
\item There exist $\Delta_B \in \C^{n , m}$, $\Delta_J \in {\rm SHerm}(n)$, $\Delta_R,\,\Delta_E \in {\rm Herm}(n)$
  satisfying~\eqref{JREBequi_1}~and~\eqref{JREBequi_2}.
\item $x$ satisfies $x_3=0$.
\end{enumerate}
Moreover, we have
\begin{eqnarray}\label{JREBequi_4}
\inf\left\{\left. {\|\Delta_J~\Delta_R ~ \Delta_E~\Delta_B\|}_F^2 \,\right|~ \Delta_J,\,\Delta_R, \, \Delta_E\in \C^{n,n},\,
\Delta_B \in \C^{n , m} ~{\rm satisfy}~\eqref{JREBequi_1}\,\text{and}~\eqref{JREBequi_2}  \right \}\nonumber
\qquad\qquad\\
=  \inf \left\{\left.\frac{{\|{\Delta}_1\|}_F^2}{2+|\lambda|^2}+{\|{\Delta}_2\|}_F^2~\right|~ {\Delta}_1\in \C^{n,n},\,
{\Delta}_2 \in \C^{n , m},\,[{\Delta}_1~{\Delta}_2] u=r,~[{\Delta}_1~{\Delta}_2]^H w=s\right \},\qquad
\end{eqnarray}
and
\begin{eqnarray}\label{JREBequi_5}
  \inf\Big\{ {\|\Delta_J~\Delta_R ~ \Delta_E~\Delta_B\|}_F^2 \,\Big|~ \Delta_J \in {\rm SHerm}(n),\,
  \Delta_R, \, \Delta_E \in {\rm Herm}(n),\, \Delta_B \in \C^{n , m} \nonumber\\
  ~{\rm satisfy}~\eqref{JREBequi_1}\,{\rm and}~\eqref{JREBequi_2}  \Big\}\nonumber \\
=  \inf \Bigg\{{\left \|\frac{{\Delta}_1 + {\Delta}_1^H}{2} \right\|}_F^2 +
  \frac{1}{1+|\lambda|^2}{\left \|\frac{{\Delta}_{1} - {\Delta}_1^H}{2} \right\|}_F^2+{\|{\Delta}_2\|}_F^2\,\Bigg |~
  {\Delta}_1\in \C^{n,n},\,{\Delta}_2 \in \C^{n , m},\nonumber\\
  \, [{\Delta}_1~{\Delta}_2] u=r,~[{\Delta}_1~{\Delta}_2]^H w=s\Bigg \}.
\end{eqnarray}
\end{lemma}

\proof ``$1) \Rightarrow 2)$'':  Let $\Delta_J$, $\Delta_R$, $\Delta_E\in \C^{n,n}$ and
$\Delta_B \in \C^{n , m}$ be such that they satisfy~\eqref{JREBequi_1} and~\eqref{JREBequi_2}.
By setting ${\Delta}_1 =\Delta_J-\Delta_R + \lambda \Delta_E$,
${\Delta}_2=\Delta_B$ and $\Delta=[{\Delta }_1~{\Delta}_2]$ we get $\Delta \in \C^{n , n+m}$
with $\Delta u=r$ and $\Delta^Hw=s$. Also, observe that by~\eqref{JREequi_6} we have
\begin{equation}
 \frac{{{\|\Delta}_1\|}_F^2}{2+|\lambda|^2} \leq {\| \Delta_J\|}_F^2+{\|\Delta_R\|}_F^2+{\| \Delta_E\|}_F^2
\end{equation}
which implies
\[
 \frac{{{\|\Delta}_1\|}_F^2}{2+|\lambda|^2}+{{\|\Delta}_2\|}_F^2
 \leq {\| \Delta_J\|}_F^2+{\|\Delta_R\|}_F^2+{\| \Delta_E\|}_F^2+{\| \Delta_B\|}_F^2.
\]
Now ``$\geq$'' in~\eqref{JREBequi_4} can be shown by arguments similar to those in the proof of
``$1) \Rightarrow 2)$'' in Lemma~\ref{thm;equivalentcond}.

\noindent
``$2) \Rightarrow 1)$'': Conversely, let $\Delta \in \C^{n , n+m}$ such that $\Delta u=r$ and ${\Delta}^Hw=s$ and
suppose that $\Delta =[{\Delta}_1~{\Delta}_2]$ where ${\Delta}_1 \in \C^{n,n}$ and $\Delta_2 \in \C^{n , m}$.
Define $$\Delta_J=\frac{\Delta_1}{2+|\lambda|^2},\quad\Delta_R=-\frac{\Delta_1}{2+|\lambda|^2},\quad
\Delta_E= \frac{\overline \lambda {\Delta}_1}{2+|\lambda|^2},\quad\mbox{and}\quad\Delta_B= {\Delta}_2,$$ then
$\Delta_J,\,\Delta_R,\, \Delta_E$ and $\Delta_B$ satisfy
$[\Delta_J -\Delta_R+\lambda \Delta_E~\Delta_B] =\Delta$ and hence~\eqref{JREBequi_1} and~\eqref{JREBequi_2}.
Furthermore, we have
\[
 {\| \Delta_J\|}_F^2+{\|\Delta_R\|}_F^2+{\| \Delta_E\|}_F^2+{\| \Delta_B\|}_F^2=
 \frac{{\|{\Delta}_1\|}_F^2}{2+|\lambda|^2}+{\|{\Delta}_2\|}_F^2.
\]
Therefore, we get ``$\leq$'' in~\eqref{JREBequi_4} by following arguments similar to those in
the proof ``$2) \Rightarrow 1)$'' in Lemma~\ref{thm;equivalentcond}.

\noindent
``$2) \Rightarrow 3)$'':
To this end, let $\Delta \in \C^{n , n+m}$ be such that $\Delta u=r$ and ${\Delta}^Hw=s$, and
suppose that $\Delta =[{\Delta}_1~{\Delta}_2]$ where ${\Delta}_1\in \C^{n,n}$ and ${\Delta }_2 \in \C^{n ,m}$.
Setting $$\Delta_R=-\frac{{\Delta}_1+{\Delta}_1^H}{2},\quad \Delta_J=\frac{{\Delta}_1-{\Delta}_1^H}
{2(1+|\lambda|^2)},\quad\Delta_E= \frac{\overline \lambda( {\Delta}_1-{\Delta}_1^H)}{2(1+|\lambda|^2)},\quad
\mbox{and}\quad \Delta_B={\Delta}_2,$$ we have $\Delta_J \in {\rm SHerm}(n)$ and, because of $\lambda\in i\mathbb R$,
also  $\Delta_R,\,\Delta_E \in {\rm Herm}(n)$. Furthermore, we obtain
\[
[\Delta_J -\Delta_R +\lambda \Delta_E\,~\Delta_B] = [{\Delta}_1 ~ {\Delta}_2]= \Delta.
\]
Thus, $\Delta_J,\Delta_R,\Delta_E$, and $\Delta_B$ satisfy~\eqref{JREBequi_1} and~\eqref{JREBequi_2}, and
we also have
\[
 {\| \Delta_J\|}_F^2+{\|\Delta_R\|}_F^2+{\| \Delta_E\|}_F^2+{\| \Delta_B\|}_F^2=
 {\left \|\frac{{\Delta}_1 + {{\Delta}_1}^H}{2} \right\|}_F^2 +
  \frac{1}{1+|\lambda|^2}{\left\|\frac{{\Delta}_1 - {{\Delta}_1}^H}{2} \right\|}_F^2
  +{\| {\Delta}_2\|}_F^2.
\]
Therefore ``$\leq$'' in~\eqref{JREBequi_5} can be shown by arguments similar to those in the
proof of ``$2)\Rightarrow 1)$'' in Lemma~\ref{lem1:str_RE}.

\noindent
``$3)\Rightarrow2)$'': Let $\Delta_R,\,\Delta_E \in {\rm Herm}(n)$, $\Delta_J \in {\rm SHerm}(n)$ and
$\Delta_B \in \C^{n , m}$ be such that they satisfy~\eqref{JREBequi_1} and~\eqref{JREBequi_2}.
Define ${\Delta}_1 =\Delta_J-\Delta_R + \lambda \Delta_E$,
${\Delta}_2=\Delta_B$ and $\Delta=[{\Delta }_1~{\Delta}_2]$ then $\Delta \in \C^{n , n+m}$
with $\Delta u=r$ and ${\Delta}^Hw=s$. Again, observe that by~\eqref{JREequi_7} we have that
\begin{equation}
{\frac{1}{1+|\lambda|^2}}{\left \|\frac{{\Delta}_1 -{{\Delta}_1}^H}{2} \right \|}_F^2
+{\left \|\frac{{\Delta}_1 +{{\Delta}_1}^H}{2} \right \|}_F^2
\leq {\|\Delta_J\|}_F^2+{\|\Delta_E\|}_F^2 +{\|\Delta_R\|}_F^2.
\end{equation}
This implies
\[
{\frac{1}{1+|\lambda|^2}}{\left \|\frac{{\Delta}_1 -{{\Delta}_1}^H}{2} \right \|}_F^2
+{\left \|\frac{{\Delta}_1 +{{\Delta}_1}^H}{2} \right \|}_F^2+ {\|{\Delta}_2\|}_F^2
\leq {\|\Delta_J\|}_F^2+{\|\Delta_E\|}_F^2 +{\|\Delta_R\|}_F^2 +{\|\Delta_B\|}_F^2,
\]
and thus, ``$\geq$'' in~\eqref{JREBequi_5} can be shown by arguments similar to those
in the proof of ``$1)\Rightarrow2)$'' in Lemma~\ref{lem1:str_RE}.

\noindent
``$2) \Leftrightarrow 4)$'': This follows immediately from Theorem~\ref{thm:spcl_matrix_map}.
\eproof

\begin{theorem} \label{thm:ustrJRBE}
Let $L(z)$ be a pencil as in~\eqref{eq:defL_cz}, and let $\lambda \in i \R$ and $x \in \C^{2n+m}\setminus \{0\}$.
Partition $x=[x_1^T~x_2^T~x_3^T]^T$ so that $x_1,\,x_2\in \C^n$ and $x_3 \in \C^{m}$ and define
$w=x_1$, $u=[x_2^T~x_3^T]^T$, $\hat u=[{(2+|\lambda|^2)^{1/2}}.x_2^T~~x_3^T]^T$,
$ r=(J-R+\lambda E)x_2 +Bx_3$, $ s=[-((J+R+\lambda E)x_1)^T~(B^Hx_1+Sx_3)^T]^T$,
and $\hat s=\big[-(2+|\lambda|^2)^{-1/2}((J+R+\lambda E)x_1)^T\,~(B^Hx_1+Sx_3)^T\big]^T$.
Then
$\eta^{\mathcal B}(J,R,E,B,\lambda,x)$ and $\eta^{\mathcal S}(J,R,E,B,\lambda,x)$
are finite if and only if $x_3=0$. Furthermore, the following statements hold.
\begin{enumerate}
 \item If $x_3=0$ then
\begin{eqnarray*}\label{eq:JREB_strberr}
 \eta^{\mathcal B}(J,R,E,B,\lambda,x) =\sqrt{ {{\|\widehat{\Delta}}_1\|}^2_F
 +\|{\widehat{\Delta}_2\|}_F^2 }
\end{eqnarray*}
and
\begin{equation*}
\eta^{\mathcal B}(J,R,E,B,\lambda)=\min\left\{\sigma_{\min}\left(\mat{cc}\frac{J-R+\lambda E}{\sqrt{2+|\lambda|^2}}& B\rix^H\right),\frac{\sigma_{\min}(J-R+\lambda E)}{\sqrt{2+|\lambda|^2}},
\right\}
\end{equation*}
where $\widehat{\Delta}_1$ and $\widehat{\Delta}_2$ are given by
\begin{equation*}\label{eq:unstr_JREBsepbrrobspecial}
[\widehat{\Delta}_1~\widehat{\Delta}_2]= \left\{\begin{array}{ll}
\frac{r\hat u^H}{\|\hat u\|^2} & \mbox{ if } x_1 = 0, \\[1ex]
\frac{w \hat s^H}{\|w\|^2} & \mbox{ if } x_2 = 0,\\[1ex]
\frac{r\hat u^H}{\|\hat u\|^2} + \frac{w\hat s^H}{\|w\|^2}\left(I_{n+m}- \frac{\hat u \hat u^H}{\|\hat u\|^2}\right)
& \mbox{ otherwise. }  \end{array}\right.
\end{equation*}
\item If $x_3=0$ then
{\small{\begin{eqnarray*}\label{eq:JREB_strberr_bound}
\sqrt{ \frac{{{\|\widetilde{\Delta}}_1\|}^2_F}{{1+|\lambda|^2}}
 +\|{\widetilde{\Delta}_2\|}_F^2 }
\leq \eta^{\mathcal S}(J,R,E,B,\lambda,x) \leq
\sqrt{{\left \|\frac{{\widetilde{\Delta}}_1+ {\widetilde{\Delta}}_1^H}{2} \right\|}_F^2 +
  \frac{1}{1+|\lambda|^2}{\left \|\frac{{\widetilde{\Delta}}_1 -
  {\widetilde{\Delta}}_1^H}{2} \right\|}_F^2+{\|{\widetilde{\Delta}}_2\|}_F^2},
\end{eqnarray*}}}
where $\widetilde{\Delta}_1$ and $\widetilde{\Delta}_2$ are given by
\begin{equation*}\label{eq:JREBsepbrrobspecial}
[\widetilde{\Delta}_1~\widetilde{\Delta}_2]= \left\{\begin{array}{ll}
\frac{ru^H}{\|u\|^2} & \mbox{ if } x_1 = 0, \\[1ex]
\frac{ws^H}{\|w\|^2} & \mbox{ if } x_2 = 0,\\[1ex]
\frac{ru^H}{\|u\|^2} + \frac{ws^H}{\|w\|^2}\left(I_{n+m} - \frac{uu^H}{\|u\|^2}\right)
& \mbox{ otherwise. }  \end{array}\right.
\end{equation*}
\end{enumerate}
\end{theorem}

\proof
The proof is analogous to that of Theorem~\ref{thm:ustrREB} by using
Lemma~\ref{lem:JREBequivalentcond} instead of Lemma~\ref{lem:REBbstr_unstr}.
\eproof

\section{Numerical experiments}

In this section, we illustrate our results with the help of numerical experiments.
In particular, we show that the backward eigenpair errors computed in the previous sections
can sometimes be significantly larger than the backward errors that correspond to perturbations
that ignore the block structure of the pencil.

In the following we compare the backward errors $\eta(L,\lambda,x)$ and
$\eta^{even}(L,\lambda,x)$ from~\eqref{17.9.17one} and~\eqref{17.9.17two} with
the block structured and symmetry structured eigenpair backward errors obtained in the
sections 4--6. We consider random pencils $L(z)$ in the form~\eqref{eq:defL_cz} and random pairs
$(\lambda,x)\in i\R \times (\C^{2n+m}\setminus\{0\})$. To make this a fair comparison it is
necessary to multiply the block structured and symmetry structured eigenpair backward errors
with a factor of $\sqrt{2}$, because each of the perturbed blocks $J$, $R$, $E$, or $B$ occurs
twice in the pencil $L(z)$.
We used Matlab Version No.~7.8.0 (R2009a) to compute the eigenpair backward errors in all cases.

\begin{example}\label{ex1:numeric}{\rm
We take a random asymptotically stable system  with $J,R,Q \in \C^{4,4}$, $B\in \C^{4,3}$,
 $S\in \C^{3,3}$ and $P=0$ such that $J^H=-J,\, R^H=R,\,Q^H=Q >0$ and $S^H=S>0$. For a particular choice of
these matrices, the corresponding pencil $L(z)$ in the form~\eqref{eq:defL_cz} turned out to have the eigenvalues
$\pm 54.518-63.914i$, $\pm 46.8738-16.2214i$, $\pm 6.8221-3.2867i$, $\pm 4.7381+11.4052i$ and $\infty$,
where $\infty$ is a semisimple eigenvalue of multiplicity $3$. Thus, the system is strictly passive. We
fix a vector $x=[x_1^T~\,x_2^T~\,x_3^T]^T\in \C^{11}\setminus\{0\}$, where $x_1,x_2 \in \C^{4}$, $0=x_3 \in \C^{3}$
and randomly select vectors $x_1$ from the intersection of the kernels of $B^H$ and $R$, and $x_2$ from the
kernel of $R$. Thus, $x$ satisfies $x_3=0$, $B^Hx_1=0$, $Rx_1=0$ and $Rx_2=0$, and hence $x$ fulfils the
finiteness criteria for all block- and symmetry-structure-preserving eigenpair backward errors from
sections 4--7.

In Table~\ref{table_1_Block}, we compare $\eta(L,\lambda,x)$ and $\eta^{even}(L,\lambda,x)$ with various
block-structure-preserving eigenpair backward errors of $L(z)$ for pairs $(\lambda,x)$, where $x$ is chosen
as above and random values for $\lambda$ on the imaginary axis were chosen. For the sake of saving space,
we omit $\lambda$ and $x$ from the notation of backward errors in Table~\ref{table_1_Block}
and also in the following Table~\ref{table_2_Sym}.

\begin{table}[!h]\small
\resizebox{.95\textwidth}{!}{\begin{minipage}{\textwidth}
\caption{\sl {Comparison of various block-structure-preserving eigenpair backward errors for the pencil $L(z)$ of Example~\ref{ex1:numeric}.}}
\begin{tabular}{|c|c|c|c|c|c|c|}
\hline
 $\lambda$ & $\eta(L)$ & $\eta^{even}(L)$ &$\sqrt{2}\,\eta^{\mathcal B}(J,E)$
 & $\sqrt{2}\,\eta^{\mathcal B}(E,B)$ & $\sqrt{2}\,\eta^{\mathcal B}(J,B)$&
 $\sqrt{2}\,\eta^{\mathcal B}(J,E,B)$ \\
 &&&$=\sqrt{2}\,\eta^{\mathcal B}(R,E)$&&$=\sqrt{2}\,\eta^{\mathcal B}(R,B)$& $=\sqrt{2}\,\eta^{\mathcal B}(R,E,B)$ \\ \hline
 .138i & 3.687 & 4.752 & 6.501 & 47.560 & 6.563 & 6.501 \\ \hline
 -.510i & 3.364 & 4.353 & 5.927 & 13.046 & 6.653 & 5.927 \\ \hline
 -.895i & 2.849 & 3.698 & 5.021 & 7.529 & 6.739 & 5.021 \\ \hline
 1.048i & 2.553 & 3.280 & 4.522 & 6.249 & 6.552 & 4.522 \\ \hline
 -1.321i & 2.346 & 3.056 & 4.139 & 5.190 & 6.859 & 4.139 \\ \hline
 1.908i & 1.734 & 2.230 & 3.095 & 3.494 & 6.668 & 3.095 \\ \hline
 2.508i & 1.405 & 1.810 & 2.524 & 2.717 & 6.817 & 2.524 \\ \hline
 \end{tabular}
\label{table_1_Block}
\end{minipage}}
\end{table}

In Table~\ref{table_2_Sym}, we record various symmetry-structure-preserving eigenpair backward errors for the same
choice of pairs $(\lambda,x)$ as in Table~\ref{table_1_Block}. We sometimes observe large differences between
various of these symmetry-structure-preserving eigenpair backward
errors. The tightness of the lower and upper bounds for $\eta^{\mathcal S}(R,E,\lambda,x)$ and
$\eta^{\mathcal S}(J,R,E,\lambda,x)$ depends on the values of $\lambda$, which is clear by
Theorem~\ref{thm:strRE} and Theorem~\ref{thm:ustrJRE}. Also
the corresponding block-structure-preserving eigenpair backward errors $\eta^{\mathcal B}(J,E,\lambda,x)$ and
$\eta^{\mathcal B}(R,E,\lambda,x)$
are sometimes significantly smaller than their symmetry-structure-preserving counterparts,
i.e., $\eta^{\mathcal S}(J,E,\lambda,x)$ and $\eta^{\mathcal S}(R,E,\lambda,x)$, respectively.

{
\begin{table}[!h]\footnotesize
\resizebox{.85\textwidth}{!}{\begin{minipage}{\textwidth}
\caption{\sl {Comparison of various symmetry-structure-pres.~eigenpair backward errors for the pencil
$L(z)$ of Example~\ref{ex1:numeric}.}}
\begin{tabular}{|c|c|c|c|c|c|c|c|}
\hline
 $\lambda$ & $\eta(L)$ & \!\!$\eta^{even}(L)$\!\! &\!\!$\sqrt{2}\,\eta^{\mathcal S}(J,E)$\!\!
 & \!{{lower bound of}}\! & \!{{upper bound of}}\! & \!{{lower bound of}}\! & \!{{upper bound of}}\!\\
 &&& & $\sqrt{2}\,\eta^{\mathcal S}(R,E)$&$\sqrt{2}\,\eta^{\mathcal S}(R,E)$&
$\sqrt{2}\,\eta^{\mathcal S}(J,R,E)$ & $\sqrt{2}\,\eta^{\mathcal S}(J,R,E)$ \\ \hline
 .138i & 3.687 & 4.752 & 8.462 & 6.563 & 38.625 & 6.501 & 6.523 \\ \hline
 -.510i & 3.364 & 4.353 & 7.647 & 6.653 & 11.330 & 5.927 & 6.178\\ \hline
 -.895i & 2.849 & 3.698 & 6.444 & 6.739 & 7.282 & 5.021 & 5.635\\ \hline
 1.048i & 2.553 & 3.280  & 5.954 & 6.249 & 6.362 & 4.522 & 5.357 \\ \hline
 -1.321i & 2.346 & 3.056 & 5.283 & 5.190 & 5.767 & 4.139 & 5.152 \\ \hline
 1.908i & 1.734 & 2.230 & 4.111 & 3.494 & 4.954 & 3.095 & 4.787 \\ \hline
 2.508i & 1.405 & 1.810 & 3.369 & 2.717 & 4.760 & 2.524 & 4.694 \\ \hline
\end{tabular}
\label{table_2_Sym}
\end{minipage}}
\end{table}
}

}
\end{example}

\section{Conclusions}


We have obtained eigenpair and eigenvalue backward errors of a pencil $L(z)$ of the form~\eqref{eq:defL_cz}
with respect to perturbations that respect the given block structure of $L(z)$ and also those that in
addition respect the symmetry structure of $L(z)$. We have shown that these backward errors may be
significantly larger than those that ignore the special block structure of the pencil.
The following table gives an overview of the existence of formulas (or bounds) for these
backward errors, when only specific blocks in the pencil are perturbed.
In the second and third column, a check mark ``$\checkmark$'' means that an explicit formula
for a block- or symmetry-structure-preserving eigenpair backward error is available for perturbations
that are restricted to blocks from the first column.
In some cases, the real eigenpair backward error is complementary. Furthermore,
in all cases block-structure-preserving eigenvalue backward errors can also be obtained while
symmetry-structure-preserving eigenvalue backward errors are obtained only for  perturbations
restricted to any two of the three blocks $J$, $R$ and $B$.

\medskip

\begin{center}
{\begin{minipage}{\textwidth}
\qquad\begin{tabular}{ |c|c|c| }
\hline
perturbed blocks & block-str.-pres. backward error& symm.-str.-pres. backward error\\ \hline
\cline{1-3}
J and E & \checkmark $\;$Theorem~\ref{thm:ustrJE}\ & \checkmark $\;$Theorem~\ref{thm:strJE} \\ \hline
R and E & \checkmark $\;$Theorem~\ref{thm:ustrRE}& bounds in Theorem~\ref{thm:strRE} \\ \hline
J and R & \checkmark $\;$Theorem~\ref{thm:ustrJR} (also real)&  \checkmark $\;$Theorem~\ref{thm:ustrJR} (also real) \\ \hline
J and B & \checkmark $\;$Theorem~\ref{thm:unstrerrorJB} (also real)& -- \\ \hline
R and B & \checkmark Remark~\ref{thm:ustrRB} (also real)& -- \\ \hline
E and B & \checkmark $\;$Theorem~\ref{thm:ustrEB} (also real)& -- \\ \hline
J,\,R and B & \checkmark $\;$Theorem~\ref{thm:ustrJRB} (also real)&\checkmark $\;$Theorem~\ref{thm:strJRB} (also real) \\ \hline
R,\,E and B & \checkmark $\;$Theorem~\ref{thm:ustrREB}& bounds in Theorem~\ref{thm:ustrREB}  \\ \hline
J,\,E and B & \checkmark $\;$Remark~\ref{thm:ustrJEB}& -- \\ \hline
J,\,R and E & \checkmark $\;$Theorem~\ref{thm:ustrJRE}& bounds in Theorem~\ref{thm:ustrJRE}\\ \hline
J,\,R,\,E and B & \checkmark $\;$Theorem~\ref{thm:ustrJRBE}&bounds in Theorem~\ref{thm:ustrJRBE}  \\ \cline{1-3} \hline
\end{tabular}
\end{minipage}}
\end{center}


\end{document}